\documentclass[letterpaper,11pt]{article}
\usepackage{amsmath}
\usepackage{amssymb}
\usepackage{setspace}
\usepackage{placeins}
\usepackage{graphicx}
\usepackage{authblk}
\usepackage{natbib}
\usepackage{lineno,mathtools}
\usepackage{multicol,multirow}
\usepackage{enumitem}
\usepackage{float}
\usepackage{lscape}
\usepackage{xcolor}
\usepackage{graphicx,wrapfig,lipsum}
\usepackage{rotating}
\usepackage[labelfont=bf, labelsep=period, width=6.5in, font=footnotesize]{caption}
\usepackage{ragged2e}
\usepackage{subcaption}
\usepackage{svg}
\usepackage[version=4]{mhchem}
\usepackage[breaklinks=true]{hyperref}
\hypersetup{
    colorlinks = true,
    citecolor={blue}
}
 
\usepackage{breakcites}
\usepackage{mathtools}  
\usepackage{xfrac}
\justifying
\DeclareMathOperator{\sech}{sech}
\usepackage{color}
\usepackage[normalem]{ulem}
\newcounter{tablefootnoteMarks}
\newcounter{tablefootnoteText}

\definecolor{darkgreen}{rgb}{0.1,.6,.1}

\usepackage{tikz}

\makeatletter
\newcommand*{\encircled}[1]{\relax\ifmmode\mathpalette\@encircled@math{#1}\else\@encircled{#1}\fi}
\newcommand*{\@encircled@math}[2]{\@encircled{$\m@th#1#2$}}
\newcommand*{\@encircled}[1]{%
  \tikz[baseline,anchor=base]{\node[draw,circle,outer sep=0pt,inner sep=.2ex] {#1};}}
\makeatother

\raggedright
\makeatletter
\renewcommand\section{\@startsection{section}{1}{0in}{-0.5\baselineskip}{0.1\baselineskip}{\normalfont\normalsize\bfseries}}
\makeatother

\makeatletter
\renewcommand\subsection{\@startsection{subsection}{1}{-\parindent}{-0.5\baselineskip}{0.1\baselineskip}{\normalfont\normalsize\textit}}
\makeatother

\setcitestyle{citesep={,},aysep={}}

\begin{document}


\begin{centering}
\Large{ 
Pace in Concert with Phase: \\Rate-induced Phase-tipping in Birhythmic Oscillators}\\
\vspace{0.5in}
Ravi Kumar K$^{1}$,  Hassan Alkhayuon$^2$, Sebastian Wieczorek$^2$ and \\ Partha Sharathi Dutta$^{1,*,}$\footnote{Email: {ravi.21maz0002@iitrpr.ac.in (Ravi Kumar K)}, {hassan.alkhayuon@ucc.ie (Hassan Alkhayuon)}, {sebastian.wieczorek@ucc.ie (Sebastian Wieczorek)}, 
$^*${Corresponding author: parthasharathi@iitrpr.ac.in (Partha Sharathi Dutta)}}\\
\vspace{0.5in}
$^1$Department of Mathematics \\ Indian Institute of Technology Ropar \\ Rupnagar, Punjab 140 001, India

\vspace{0.15in}
$^2$School of Mathematical Sciences\\ University College Cork\\
Cork T12 XF62, Ireland\\
\vspace{0.4in}
\end{centering}

\doublespacing
\newpage
\justifying



\section*{Abstract}

We study rate-induced phase-tipping (RP-tipping) between two stable limit cycles of a birhythmic oscillator. We say that such an oscillator RP-tips when a time variation of an input parameter preserves the bistability of the limit cycles but induces transitions from one stable limit cycle to the other, causing abrupt changes in the amplitude and frequency of the oscillations. Crucially, these transitions occur when: the rate of change of the input is in a certain interval bounded by {\em critical rate(s)}, and the system is in {\em certain phases} of the cycle. 

We focus on two illustrative examples: the birhythmic van der Pol oscillator and the birhythmic Decroly-Goldbeter glycolysis model, each subjected to monotone and \linebreak non-monotone shifts in their input parameters. We explain RP-tipping in terms of properties of the autonomous frozen system, including the {\em phase} of a cycle and {\em partial basin instability} along the parameter path traced by the changing input. We show that RP-tipping can occur as an irreversible one-way transition or as a series of transitions between the stable limit cycles. Finally, we present RP-tipping diagrams showing combinations of the rate and magnitude of parameter shifts and the phase of the oscillation that give rise to this genuine non-autonomous instability.\\

\noindent \textbf{Keywords:} RP-tipping, birhythmicity, time-dependent external input, non-autonomous systems, phase of a limit cycle, partial basin instability

\newpage

\section{Introduction}

Critical transitions are nonlinear phenomena that can be understood as abrupt, significant, and unforeseen changes in the state of a dynamical system, often accompanied by crossing a tipping point due to minor and gradual changes in the external input \citep{scheffer2009critical,ashwin2012tipping}. The notion of a tipping point was popularised by \citet{gladwell2000tipping} and has since found applications in diverse fields like ecology \citep{scheffer2001catastrophic,carpenter2011early,dakos2014critical,o2023early}, climate science \citep{ashwin2012tipping,ritchie2021overshooting,lenton2008tipping}, and systems biology \citep{trefois2015critical,sarkar2019anticipating}. Given the potential for sudden and often irreversible changes in natural and human systems linked to these tipping points, it becomes crucial to identify and comprehend the underlying dynamical mechanisms that trigger such transitions. Based on the trigger mechanism, tipping points can be categorised into four types \citep{ashwin2012tipping,halekotte2020}: bifurcation-induced tipping (B-tipping),  rate-induced tipping (R-tipping), noise-induced tipping (N-tipping), and shock-induced tipping (S-tipping)\footnote{In some sense, one can view N-tipping and S-tipping as special cases of R-tipping.}.  Furthermore, depending on the form of the {\em base state} from which the critical transition occurs, that is, the attractor for the {\em frozen system} with fixed-in-time external input, one can speak of partial-tipping or phase-tipping (P-tipping) from limit cycles \citep{alkhayuon2018rate,alkhayuon2021phase}, multi-frequency tipping (M-tipping) from quasiperiodic torus \citep{keane20}, and fragmented-tipping (F-tipping) from spatial patterns \citep{bastiaansen2022fragmented}.
In this paper, changing {\em external inputs} are represented by shifting {\em input parameters} in the oscillator models, and are sometimes referred to as {\em external forcing}.

In general, R-tipping occurs when the input parameter changes at rates that are both fast enough and optimal, so the system fails to track the {\em moving base state} (quasi-static attractor), crosses some {\em moving threshold} in the phase space, and transitions to an {\em alternative state}. In contrast to B-tipping, R-tipping can (and usually does) occur even though the changing input parameter does not cross any autonomous bifurcation points of the base state in a frozen system. Optimal rates mean that R-tipping may occur above some critical rate or between two critical rates of parameter change \citep{o2020tipping,longo2021,ritchie2023rate}. The moving threshold can be a moving basin boundary \citep{o2020tipping} or a moving excitability quasi-threshold \citep{o2023rate}. The alternative state can be long-lived (another attractor in a multi-stable system) \citep{vanselow2022evolutionary,chaparro2021fast}, or finite-time (a transient response in an excitable system) \citep{o2023rate,kaur2022critical,vanselow2024rate}. The identification of critical rates and the exploration of the outcomes when these rates are crossed hold significant importance in the study of natural and human systems \citep{wieczorek2011excitability,van2021understanding,kaur2022critical,o2023rate,ritchie2023rate}. Most studies on R-tipping have considered sudden transitions from a stationary base state (a stable equilibrium) \citep{wieczorek2023rate,ritchie2023rate}. However, many complex systems \citep{mickens1996oscillations,goldbeter1997biochemical} exhibit non-stationary base states such as stable limit cycles, tori, or even chaotic attractors.  Sudden transitions from non-stationary base states have been less studied and understood \citep{alkhayuon2018rate,alkhayuon2021phase,ashwin2021physical,longo2021,duenas2023critical, longo2024critical}.

P-tipping occurs in nonlinear systems with an oscillatory base state, when an external input causes the system to transition to an alternative state, but only from certain phases of the oscillations \citep{alkhayuon2018rate,alkhayuon2021phase,longo2024critical}. P-tipping stands for either partial tipping or phase tipping. \citet{alkhayuon2018rate} investigated the occurrence of R-tipping from a stable limit cycle using the framework of pullback attractors in non-autonomous systems with asymptotic parameter shifts. They provided an illustrative example where the pullback attractor of the system fails to track the moving limit cycle base state (quasistatic attractor) of the frozen system partially. In other words, some trajectories of the pullback attractor diverge while others track the moving limit cycle base state. Crucially, this is a typical scenario for loss of tracking of moving non-stationary base states. Thus, this phenomenon was called {\em partial tipping} in~\citet{alkhayuon2018rate}. In an alternative approach, Longo et al. \citep{longo2021,longo2024critical} and Due\~nas et al. \citep{duenas2023critical} studied the concept of partial tipping on the hull of a non-autonomous scalar ordinary differential equation. Recently, \citet{alkhayuon2021phase} studied tipping from a stable predator-prey limit cycle to an extinction equilibrium in two paradigmatic non-autonomous predator-prey models with an Allee effect and a time-varying input in the form of (discontinuous) jumps in the growth rate of the prey. In this setting, the tipping depends purely on the phase of the cycle when the jump occurs. Hence, it was called {\em phase-tipping}.

RP-tipping is a combination of R-tipping and P-tipping in the sense that the system tips from an oscillatory base state if two conditions are satisfied simultaneously: 
(i) the rate of change of the input is both fast enough and optimal, and
(ii) the system state is in certain phases of the oscillatory base state.

One novelty of our work is that we explore RP-tipping from a base limit cycle to an alternative stable limit cycle. Such tipping implies a shift in the oscillation characteristics, e.g., amplitude and frequency. For example, an energy harvesting system exhibits the coexistence of two stable limit cycles in its wind-induced vibration \citep{kwuimy2015recurrence,zhang2023bifurcations}. There, an undesired transition to the stable limit cycle with a smaller amplitude and different frequency results in less mechanical deformation and, consequently, less efficient conversion from wind to electrical energy. Here, we illustrate RP-tipping in two birhythmic oscillators: the birhythmic van der Pol oscillator \citep{kaiser1991bifurcation}, which incorporates a higher-order polynomial damping factor, and the Decroly-Goldbeter glycolysis model of biochemical oscillations \citep{moran1984onset}. Both systems exhibit a coexistence of two stable limit cycles separated by an unstable limit cycle, which is the defining feature of birhythmic oscillators \citep{goldbeter1997biochemical,biswas2016control}.

To study RP-tipping in the aforementioned birhythmic oscillators, we subject these oscillators to monotone and non-monotone continuous time-dependent external inputs. Specifically, we consider parameter shifts \citep{ashwin2017parameter,wieczorek2021compactification,ritchie2023rate} that do not cross any autonomous bifurcation points of the frozen system. In other words, there are always two stable limit cycles at any fixed level of the input parameter. The transitions between stable limit cycle attractors happen from certain phases, known as tipping phases. Hence, it is important to define the phase of the limit cycle explicitly. It is also essential to note that the inclusion of time-dependent external inputs makes the systems non-autonomous, meaning that the traditional autonomous bifurcation theory proves inadequate for exploring RP-tipping between stable limit cycles. Recent decades have witnessed the development of various methods \citep{kloeden2011nonautonomous,ghil2023dynamical,wieczorek2021compactification,wieczorek2023rate,anagnostopoulou2023nonautonomous} to study the dynamics of non-autonomous systems. One among these methods is "basin instability," which uses the geometric properties of the frozen autonomous system \citep{o2020tipping,ritchie2023rate}. In \citet{wieczorek2023rate}, the authors showed that the notion of basin instability provides easily verifiable sufficient conditions for the occurrence of R-tipping from equilibrium base states in higher dimensional nonlinear systems. \citet{alkhayuon2021phase} extended this notion to "partial basin instability" to study P-tipping from oscillatory base states. In their work, partial basin instability of a stable limit cycle occurs by crossing the stable manifold of a saddle equilibrium, whereas, in our case, it happens by crossing the repelling limit cycle. We determine the region of partial basin instability of the coexisting stable limit cycles in their corresponding autonomous frozen system, and we identify appropriate parameter shifts that cause RP-tipping. In particular, we analyse two types of RP-tipping: one-way transitions and a series of alternating transitions between the two stable limit cycles.

We start with one-way transitions that can be interpreted as irreversible RP-tipping for both monotone and symmetric non-monotone parameter shifts. These transitions occur from the larger-amplitude base limit cycle to the smaller-amplitude alternative stable limit cycle in the birhythmic van der Pol oscillator, and conversely in the glycolysis model. We obtain two-dimensional RP-tipping diagrams with variations in the magnitude and rate of the parameter shifts. Depending on the model and the type of parameter shift, we observe either bounded or half-bounded RP-tipping regions, leading to the identification of two critical rates or one critical rate, respectively, for a given magnitude of the shift.

Next, we demonstrate a series of alternating transitions between both stable limit cycles in the birhythmic van der Pol oscillator with a non-monotone parameter shift in the form of an impulse. In this sequence, the system also transitions from the smaller-amplitude stable limit cycle to the larger-amplitude limit cycle and then back to the smaller-amplitude limit cycle.

The paper is organised as follows. Section~\ref{Sec: MathMod} introduces two autonomous birhythmic oscillator models and identifies the coexistence of two stable limit cycles in each model using one- and two-parameter bifurcation diagrams, phase portraits, and the corresponding time series. Then, we set up the non-autonomous problem to study RP-tipping (subsection~\ref{SSec: time-varying input}). In section~\ref{Sec: phase_and_BI}, we define the phase of limit cycles and relate the concept of basin instability to the RP-tipping phenomenon. In section~\ref{Sec: Result}, we demonstrate RP-tipping in the birhythmic van der Pol model and the Decroly-Goldbeter glycolysis model. Finally, we summarise and conclude in section~\ref{Sec: conclusion}.

\section{Mathematical models of birhythmic oscillators \label{Sec: MathMod}}

In this section, we introduce the two distinct birhythmic oscillator models, identify the time-varying input parameters, and highlight limit cycle bistability in the corresponding frozen systems with fixed-in-time input parameters.

\subsection{Model 1: The birhythmic van der Pol model \label{SSec: birhythmic vdp}}

\begin{figure}[ht!]
\centering
\includegraphics[width=0.9\textwidth]{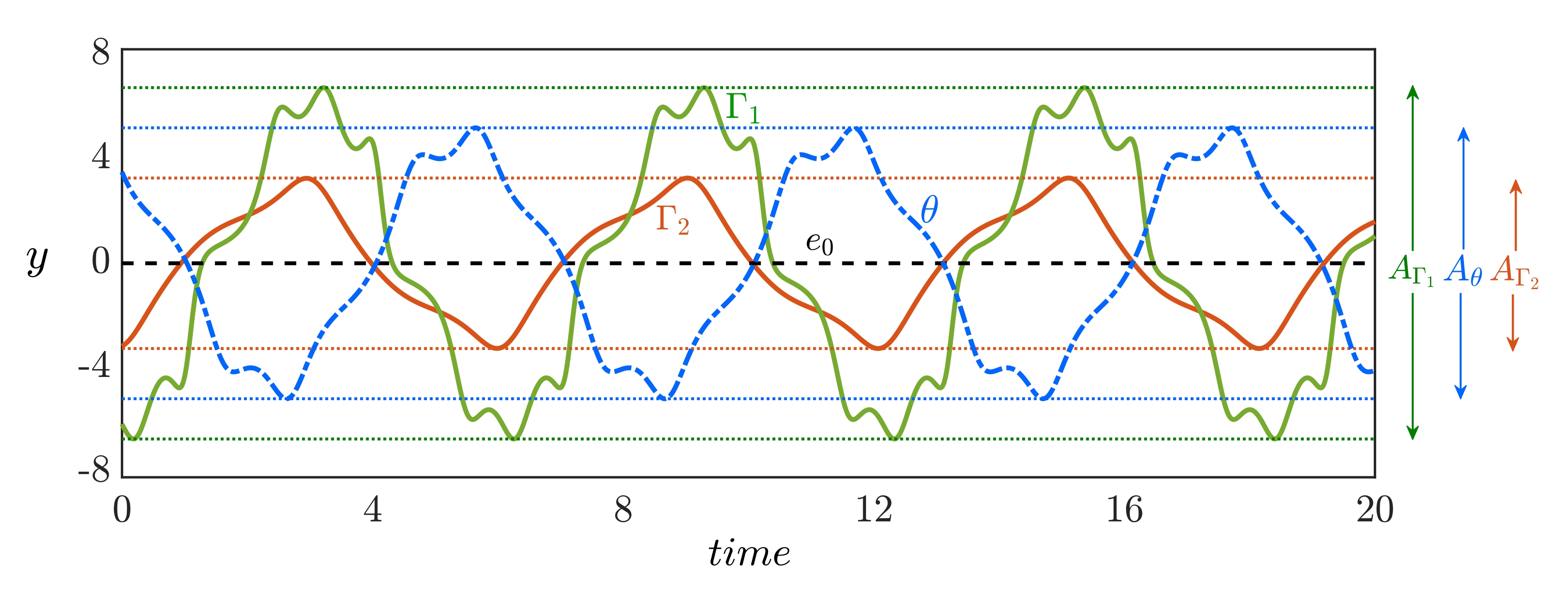}
\caption{Time series of the birhythmic van der Pol oscillator \eqref{vdp_with_xy} exhibiting limit cycle oscillations with different amplitudes and frequencies. The horizontal dotted lines indicate the maxima and minima of the large-amplitude stable limit cycle $\Gamma_1$ (green), the unstable limit cycle $\theta$ (blue), and the small-amplitude stable limit cycle $\Gamma_2$ (red) along with the amplitudes $A_{\Gamma_1}$, $A_{\theta}$, and $A_{\Gamma_2}$ (marked with double sided arrows), respectively. The dashed black line at $y = 0$ denotes an unstable equilibrium point $e_0$. The parameter values are $d = -0.1, \alpha = 0.114$, $\beta = 0.003$, and $\mu = 0.6$.} \label{fig: Figure_timeseries_vdp}
\end{figure}

\begin{figure}[ht!]
\centering
\includegraphics[width=0.9\textwidth]{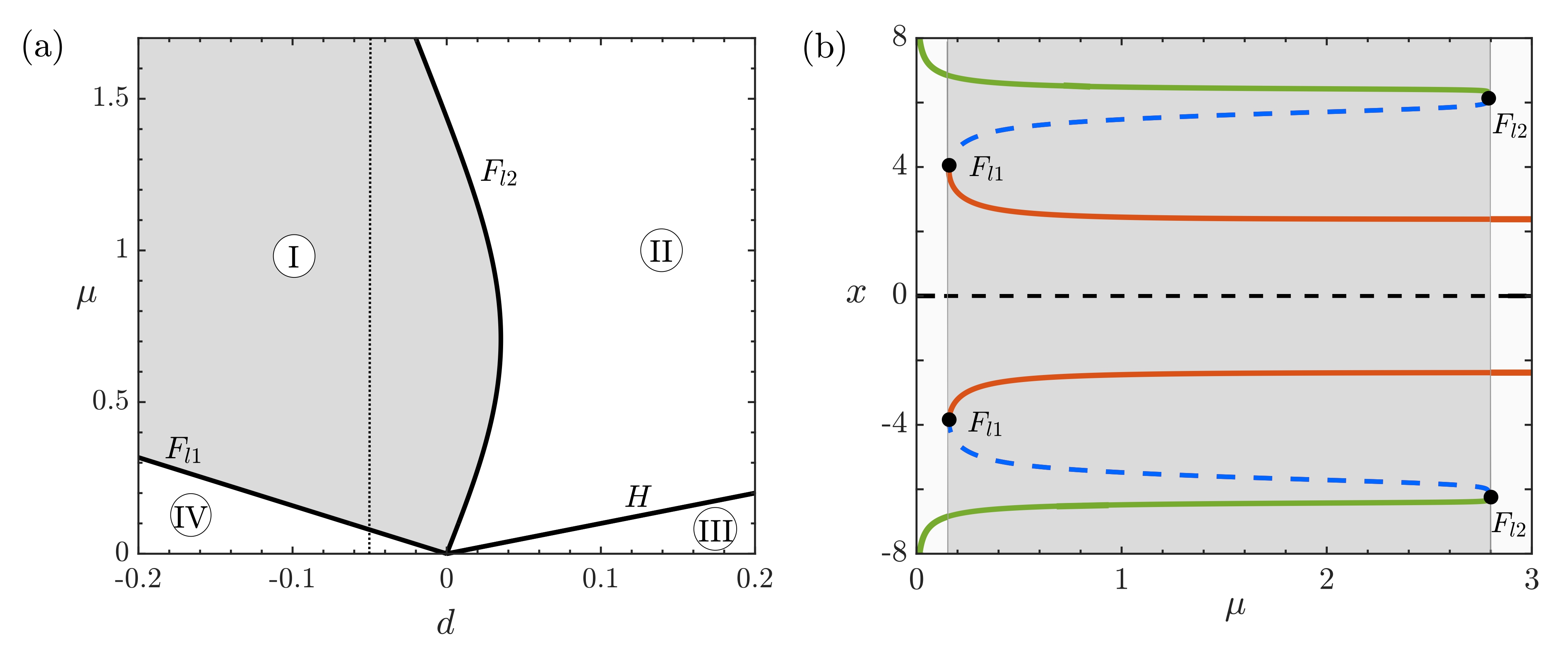}
\caption{(a) Two-parameter bifurcation diagram of the birhythmic van der Pol oscillator \eqref{vdp_with_xy} with variations in the parameters $d$ and $\mu$. $F_{l1}$, and $F_{l2}$ represent the fold/saddle-node bifurcation of limit cycles (SNLC) curves where a stable limit cycle and an unstable limit cycle merge and disappear, and $H$ represents the Hopf bifurcation curve. The grey-shaded region-I marks the coexistence of two stable limit cycles $\Gamma_1$ and $\Gamma_2$ (birhythmicity) separated by an unstable limit cycle $\theta$; in region-II, the small amplitude stable limit cycle $\Gamma_2$ exists alone (monorhymicity); region-III marks the existence of a stable steady state (point attractor); and in region-IV, the large amplitude stable limit cycle $\Gamma_1$ exists alone (monorhymicity). The vertical dotted line represents the value of $d=-0.05$ for which the one-parameter bifurcation diagram is plotted. (b) One-parameter bifurcation diagram exhibiting the maxima and minima of each limit cycle with variations in the parameter $\mu$. The stable limit cycles are marked with solid (red and green) curves, and the unstable limit cycle is marked with a dashed (blue) curve. The dashed horizontal line at $x=0$ represents an unstable equilibrium point. The grey-shaded region represents the existence of birhythmicity. Other parameter values are $\alpha = 0.093$, and $\beta = 0.0019$.}
\label{fig: Figure_bifn_vdp}
\end{figure}

The birhythmic van der Pol oscillator equation for a real variable $x$,
\begin{linenomath*}
\begin{equation}
     \ddot{x} -\mu(1-x^2+\alpha x^4-\beta x^6)\dot{x}+x + d(\dot{x} - x)=0,
     \label{vdp_with_d}
\end{equation}
\end{linenomath*}
was introduced by Kaiser \citep{kaiser1991bifurcation}, who modified the classical van der Pol oscillator equation with an additional higher-order polynomial damping factor $\mu(\alpha x^4-\beta x^6)\dot{x}$. Later, \citet{biswas2016control} included a new feedback term $d(\dot{x} - x)$. For our analysis, it is convenient to rewrite the second-order differential equation \eqref{vdp_with_d} as two coupled first-order differential equations: 
\begin{linenomath*}
\begin{equation}
    \label{vdp_with_xy}
    \begin{cases} 
        \dot{x} &= y,\\ 
        \dot{y} &= \mu(1-x^2+\alpha x^4-\beta x^6)y-x-d(y-x),\\
     \end{cases}
\end{equation}
\end{linenomath*}
where $\mu > 0$ is the {\em external input parameter}, $\alpha > 0$ and $\beta > 0$ are non-linear damping parameters, and $d$ is the feedback strength. In Ref. \citep{biswas2016control}, the authors assumed 
that the periodic solutions of \eqref{vdp_with_d} can be approximated by
$x(t) = A\cos(\omega t)$, 
where $A$ is the amplitude and $\omega$ is the angular frequency of the oscillator, and then by employing the harmonic decomposition method they derived the following
amplitude equation:
\begin{linenomath*}
\begin{equation}
\mu\left(1-\dfrac{1}{4}A^2+\dfrac{\alpha}{8}A^4-\dfrac{5\beta}{64}A^6\right) - d= 0.
\label{amp_eqn_w_o_d}
\end{equation}
\end{linenomath*}
Equation \eqref{amp_eqn_w_o_d} provides an approximate relationship between the system parameters and the oscillation amplitude. The number of roots of \eqref{amp_eqn_w_o_d} for a given $\mu$ and $d$ determines the number of possible periodic solutions to Eq.~\eqref{vdp_with_d}, (see Fig.~\ref{fig: Figure_timeseries_vdp} and  Fig.~\ref{fig: Figure_bifn_vdp}(b)).

There is a range of parameters where \eqref{amp_eqn_w_o_d} has three positive real roots. Two of these roots, namely $A_{\Gamma_{1}} > A_{\Gamma_{2}}$,  correspond to two stable limit cycles $\Gamma_{1}$ and $\Gamma_{2}$. The third root, namely $A_{\theta}$, where $A_{\Gamma_{1}} > A_{\theta} > A_{\Gamma_{2}}$, corresponds to the unstable limit cycle $\theta$. The time series for all three limit cycles are shown in Fig.~\ref{fig: Figure_timeseries_vdp}. 

To identify the range of parameters where all three limit cycles exist, we fix $\alpha=0.093$ and $\beta=0.0019$ and produce a two-parameter bifurcation diagram for the autonomous model \eqref{vdp_with_xy} in the plane $(d,\mu)$ in Fig.~\ref{fig: Figure_bifn_vdp}(a). The bifurcation analyses have been performed using the \textit{XPPAUT} package \citep{xppaut}.  In Fig.~\ref{fig: Figure_bifn_vdp}(a), the grey-shaded region-I marks the parameter region with three limit cycles. The maxima and minima for two stable and one unstable limit cycles are shown in the grey-shaded portion of the one-parameter bifurcation diagram in Fig.~\ref{fig: Figure_bifn_vdp}(b) for $d=-0.05$ (corresponding to the vertical dashed line in Fig.~\ref{fig: Figure_bifn_vdp}(a)). Depending on the chosen initial conditions, the system's trajectories converge towards either the large-amplitude ($\Gamma_1$) or the small-amplitude ($\Gamma_2$) stable limit cycle, represented by the green and red branches, respectively, in Fig.~\ref{fig: Figure_bifn_vdp}(b). The blue dashed branches represent the maxima and minima of the unstable limit cycle ($\theta$), acting as a separatrix between the two stable limit cycles.  Pairs of limit cycles collide and disappear at the fold (saddle-node) of limit cycles bifurcation points $F_{l1}$ and $F_{l2}$. The dashed black line at $x=0$ represents the unstable equilibrium $e_0$.

\begin{figure}[ht!]
\centering
\includegraphics[width=0.925\textwidth]{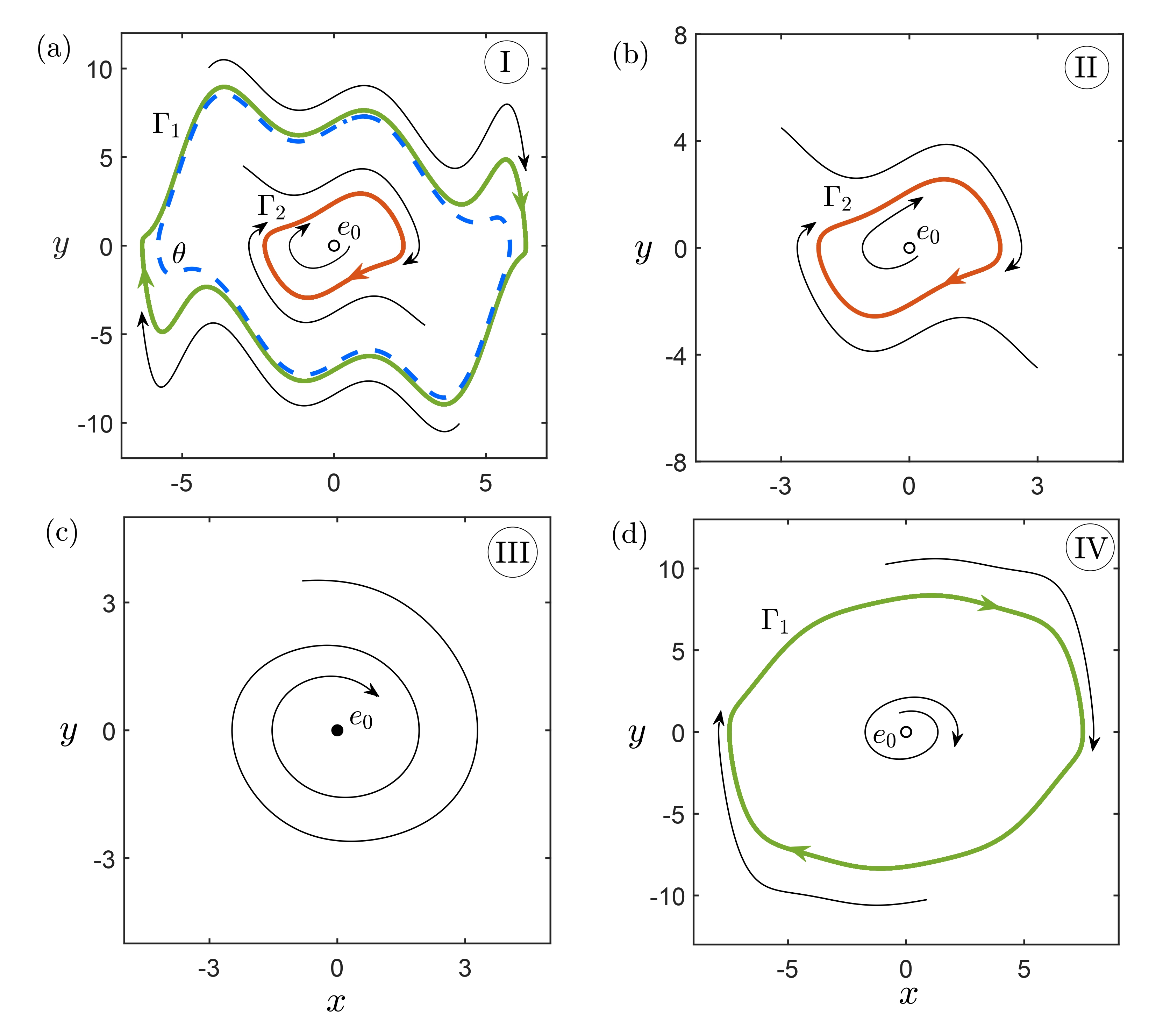}
\caption{Phase portraits of the birhythmic van der Pol oscillator \eqref{vdp_with_xy} corresponding to distinct regions shown in Fig.~\ref{fig: Figure_bifn_vdp}(a): For (a) $d = -0.001$ and $\mu = 1$ (region I: birhythmicity); (b) $d = 0.1$ and $\mu = 1$ (region II: monorhythmicity - small amplitude stable limit cycle $\Gamma_2$); (c) $d = 0.15$ and $\mu = 0.05$ (region III: stable equilibrium point); and (d) $d = -0.15$ and $\mu = 0.05$ (region IV: monorhythmicity - large amplitude stable limit cycle $\Gamma_1$). $e_0$ marked with an open circle stands for an unstable equilibrium point, and a filled circle stands for a stable equilibrium point. The blue dashed curve $\theta$ stands for the unstable limit cycle. The black curves headed with an arrow represent some exemplary trajectories. Other parameter values are same as in Fig.~\ref{fig: Figure_bifn_vdp}.}  \label{fig: Figure_phaseport_vdp}
\end{figure}

To each side of the birhythmic region-I in Fig.~\ref{fig: Figure_bifn_vdp}(a), there is a monorhythmic region with a single small-amplitude stable limit cycle $\Gamma_2$ in region-II, or single large-amplitude stable limit cycle $\Gamma_1$ in region-IV. The small-amplitude stable limit cycle $\Gamma_2$  from region-II disappears along the supercritical Hopf bifurcation curve $H$, leaving just one stable equilibrium point $e_0$ in region-III. Examples of phase portraits in different regions from Fig.~\ref{fig: Figure_bifn_vdp}(a), including the cycles $\Gamma_1$, $\Gamma_2$, $\theta$, and the equilibrium $e_0$ (stable in region-III and unstable in other regions), are shown in Fig.~\ref{fig: Figure_phaseport_vdp}.

\begin{figure}[ht!]
\centering
\includegraphics[width=0.9\textwidth]{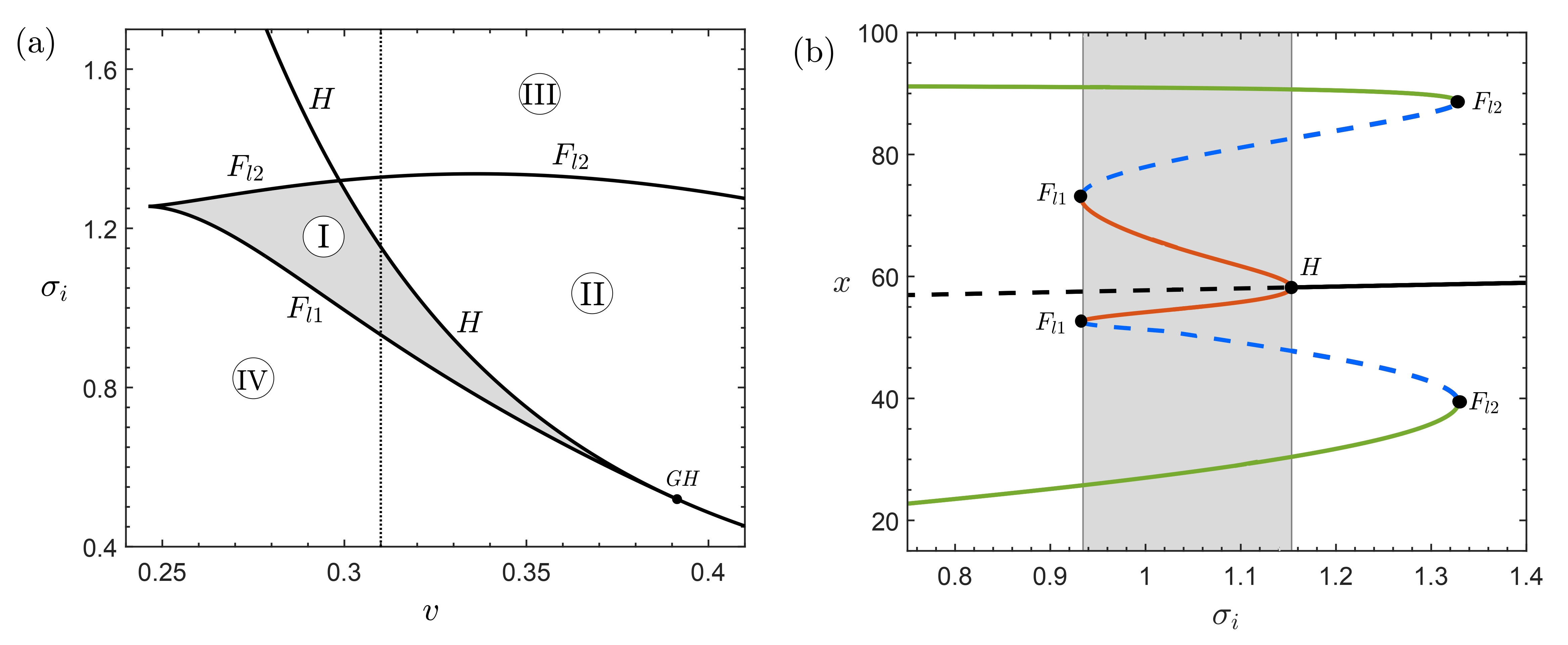}
\caption{(a) Two-parameter bifurcation diagram of the Decroly-Goldbeter glycolysis model \eqref{gly_model} with variations in the parameters $v$ and $\sigma_i$. $F_{l1}$, and $F_{l2}$ represent the SNLC curves where a stable limit cycle and an unstable limit cycle merge and disappear, and $H$ represents the Hopf bifurcation curve. The filled circle marked by GH at $(v,\sigma_i) \approx (0.3914,0.5195)$ represents the generalised Hopf (Bautin) bifurcation point where $F_{l1}$ and $H$ merge. The grey-shaded region-I marks the coexistence of two stable limit cycles $\Gamma_1$ and $\Gamma_2$, separated by an unstable limit cycle $\theta$; in region-II, the large amplitude stable limit cycle $\Gamma_1$ coexists with a stable equilibrium point (hard-excitation); region-III marks the existence of a stable equilibrium point; and in region-IV, a single stable limit cycle exists. The vertical dotted line represents $v=0.31$, for which the one-parameter bifurcation diagram is plotted. (b) One-parameter bifurcation diagram with variations in the input parameter $\sigma_i$. The maxima and minima of stable limit cycles are marked with solid (red and green) curves, and those of the unstable limit cycle are marked with dashed (blue) curves. The almost horizontal dashed (black) line represents an unstable equilibrium point, and the solid (black) line represents a stable equilibrium point. The grey-shaded region represents the existence of birhythmicity. Other parameter values are $K = 10$, $L = 3.6 \times 10^6$, $\sigma_M = 10$, $n = 5$, $q = 1$, and $k_s = 0.06$.}
\label{fig: Figure_bifn_gly}
\end{figure}

\subsection{Model 2: The Decroly-Goldbeter glycolysis model \label{SSec: glycolysis model}}

The Decroly-Goldbeter glycolysis model of enzyme reaction \citep{moran1984onset} in two real variables $x$ (substrate) and $y$ (product) is given by the following equations:
\begin{linenomath*}
\begin{equation}
\label{gly_model}
\begin{cases} 
\dot{x} &= v + \dfrac{\sigma_{i}}{K^n+y^n}\,y^n-\sigma_{M}\,\Phi(x,y),\\ 
\dot{y} &= q\,\sigma_{M}\,\Phi(x,y) - k_{s}\,y - q\,\dfrac{\sigma_{i}}{K^n+y^n}\,y^n.
\end{cases}
\end{equation}
\end{linenomath*}
Here, the parameter $v$ is the constant input of substrate $x$, the {\em input parameter} $\sigma_i$ is the normalized maximum rate of the recycling enzymes,  $\sigma_M$ is the normalized maximum rate of reaction, $q = M/K_p$, where $M$ is the Michaelis constant and $K_p$ denotes the dissociation constant of the product, and $k_s$ is the removal rate of the product $y$. The term $\dfrac{\sigma_{i}}{K^n+y^n}\,y^n$ describes the feedback mechanism of the product into the substrate with the Hill coefficient $n \geq 3$ and the half-saturation constant of the product $K$.  The rate of reaction function is given as \citep{monod1965nature}:
\begin{linenomath*}
\begin{equation}
\label{rate_of_reaction}
\Phi(x,y) = \dfrac{x(1+x)(1+y)^2}{L+(1+x)^2(1+y)^2}, \nonumber
\end{equation}
\end{linenomath*}
where $L$ denotes the allosteric constant of the enzyme.

\begin{figure}[ht!]
\centering
\includegraphics[width=0.9\textwidth]{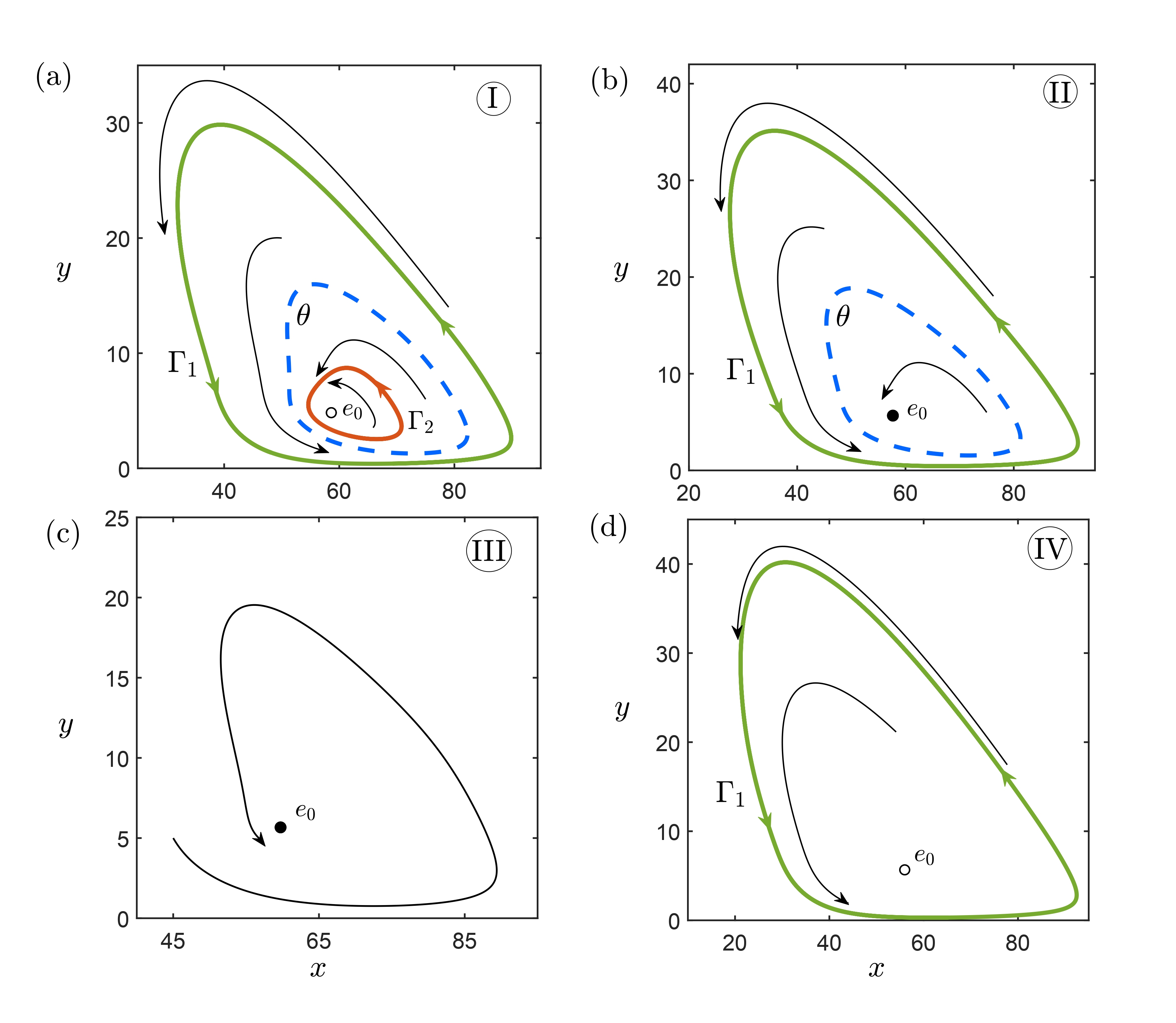}
\caption{Phase portraits illustrating different dynamics of the Decroly-Goldbeter glycolysis model~\eqref{gly_model}, corresponding to distinct regions as outlined in Fig.~\ref{fig: Figure_bifn_gly}(a): For (a) $v = 0.29$ and $\sigma_i = 1.15$ (region I: birhythmicity); (b) $v = 0.34$ and $\sigma_i = 1.1$ (region II: hard excitation); (c) $v = 0.34$ and $\sigma_i = 1.6$ (region III: stable equilibrium point); and (d) $v = 0.34$ and $\sigma_i = 0.7$ (region IV: monorhythmicity). $e_0$ marked with an open circle stands for an unstable equilibrium point, and a filled circle stands for a stable equilibrium point. The blue dashed curve $\theta$ stands for the unstable limit cycle. The black curves headed with an arrow represent some exemplary trajectories. Other parameters are same as in Fig.~\ref{fig: Figure_bifn_gly}.}
\label{fig: Figure_phaseport_gly}
\end{figure}

Fig.~\ref{fig: Figure_bifn_gly} shows the bifurcation structures of the autonomous model \eqref{gly_model}. The two-parameter bifurcation diagram (Fig.~\ref{fig: Figure_bifn_gly}(a)) with variations in the parameters $v$ and $\sigma_i$ identifies four qualitatively different regions, where the grey-shaded region-I marks the coexistence of two stable and one unstable limit cycles.  We denote the stable limit cycles with $\Gamma_{1}$ and $\Gamma_{2}$, the unstable limit cycle with $\theta$, and use the same colour coding as for the birhythmic van der Pol oscillator. Region-II presents hard excitation representing the coexistence of one stable limit cycle, one unstable limit cycle, and a stable equilibrium point, region-III marks the presence of a single stable equilibrium point, and region-IV denotes the existence of a single stable limit cycle. The region-I is bounded by one Hopf bifurcation curve $H$ and two saddle-node bifurcations of limit cycle curves $F_{l1}$ and $F_{l2}$. The maxima and minima for two stable and one unstable limit cycles are shown in the one-parameter bifurcation diagram in Fig.~\ref{fig: Figure_bifn_gly}(b) for $v=0.31$ (corresponding to the vertical dashed line in Fig.~\ref{fig: Figure_bifn_gly}(a)). The grey-shaded region in Fig.~\ref{fig: Figure_bifn_gly}(b) marks the coexistence of two stable and one unstable limit cycles.

Examples of phase portraits for different regions from Fig.~\ref{fig: Figure_bifn_gly}(a), including the cycles $\Gamma_1$, $\Gamma_2$, $\theta$, and the equilibrium $e_0$, are shown in Fig.~\ref{fig: Figure_phaseport_gly}.


\subsection{ Time-varying external inputs and non-autonomous models \label{SSec: time-varying input}}

In the remainder of the paper, we will analyse non-autonomous versions of models 
\eqref{vdp_with_xy} and \eqref{gly_model}, where the external input parameters, $\mu$ and $\sigma_i$, are replaced with time-varying external inputs, 
\begin{linenomath*}
$$\mu(rt)\quad\mbox{and}\quad \sigma_i(rt),$$ 
\end{linenomath*}
respectively. Here, the {\em rate parameter} $r > 0$ quantifies the timescale of the varying external inputs. We will be interested in how the response of each birhythmic oscillator to its external input changes as $r$ is varied.

\begin{figure}[ht!]
\centering
\includegraphics[width=0.95\textwidth]{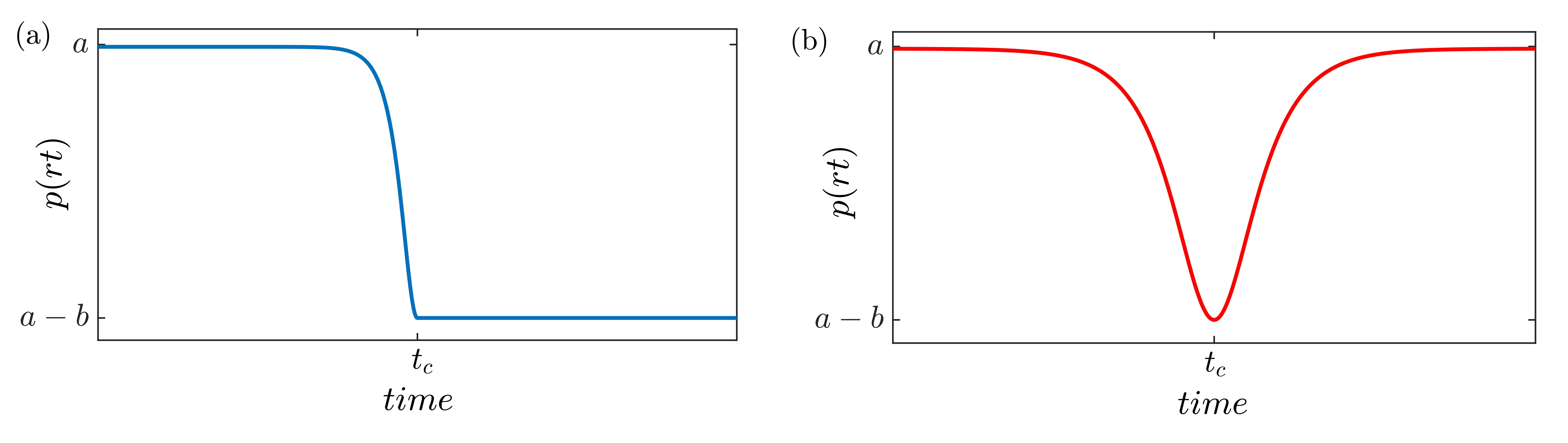}
\caption{Examples of (a) monotone and (b) non-monotone time-dependent external inputs \eqref{monotone_shift} and \eqref{nonmonotone_shift}, respectively.
} \label{fig: Figure_parameter_shape}
\end{figure}

To simplify the discussion, we use $p(rt)$ to denote either $\mu(rt)$ or $\sigma_i(rt)$, and $X$ to denote $(x,y)$ from models \eqref{vdp_with_xy} and \eqref{gly_model}. In other words, now each model can be written as:
\begin{linenomath*}
\begin{equation}
    \label{eq:nonauto_pbm_setting}
    \dot{X} = \dfrac{dX}{dt} = g(X,p(rt)).
\end{equation}
\end{linenomath*}
Furthermore, we will consider bi-asymptotically constant external inputs~\citep{wieczorek2023rate}, meaning that,
\begin{linenomath*}
\begin{equation}
    \label{limits}
    \displaystyle \lim_{t \to +\infty} p(rt) = p^+ \in \mathbb{R}, ~~ \text{and} ~~ 
    \displaystyle \lim_{t \to -\infty} p(rt) = p^- \in \mathbb{R}.
\end{equation}
\end{linenomath*}
To be specific, we consider monotone external input of the form:
\begin{linenomath*}
\begin{equation}
    \label{monotone_shift}
    p(rt) = \begin{cases} a - b\sech (r(t-t_c)) &~\text{if}~~~ t \leq t_c,\\
    a - b &~\text{if} ~~~ t > t_c,\end{cases}
\end{equation}
\end{linenomath*}
and the corresponding non-monotone external input of the form:
\begin{linenomath*}
\begin{equation}
    \label{nonmonotone_shift} 
    p(rt) = a - b\sech (r(t-t_c)), 
\end{equation}
\end{linenomath*}
shown in Figs.~\ref{fig: Figure_parameter_shape}(a) and \ref{fig: Figure_parameter_shape}(b), respectively.

Furthermore, we refer to system~\eqref{eq:nonauto_pbm_setting} with a fixed in-time input parameter $p(rt) = p$, 
\begin{linenomath*}
\begin{equation}
    \label{eq:auto_pbm_setting}
    \dot{X} = \dfrac{dX}{dt} = g(X,p),
\end{equation}
\end{linenomath*}
as the {\em autonomous frozen system}. Note that when $p = p^- (\mbox{or}\; p^+)$, then the autonomous frozen system is referred to as the past (or future) limit system. Finally, we write
\begin{linenomath*}
$$
X(t,t_0,X_0),
$$
\end{linenomath*}
to denote the solution of system~\eqref{eq:nonauto_pbm_setting} at time $t$ started from $X_0$ at initial time $t_0$. We also write 
\begin{linenomath*}
$$
X(t,X_0;p),
$$
\end{linenomath*}
to denote the solution of system~\eqref{eq:auto_pbm_setting} at time $t$ started from initial condition $X_0$ at time $t=0$ for a fixed-in-time input parameter $p$.

\section{When to expect rate-induced phase-tipping \label{Sec: phase_and_BI}}

In this section, we discuss the concept of partial basin instability for limit cycles \citep{alkhayuon2021phase} and show that it is a useful criterion for testing systems’ susceptibility to RP-tipping. We start by characterising all points on a limit cycle by phases.

\subsection{Phase of limit cycles \label{SSec: PhaseLC}}
The solution $X(t,\gamma_0;p)$ to the autonomous frozen system~\eqref{eq:auto_pbm_setting} started from some $X_0 = \gamma_0 \in \mathbb{R}^2$ for a fixed $p$ is periodic if there exists a constant $T>0$ such that 
\begin{linenomath*}
$$
X(t,\gamma_0;p) = 
X(t + T,\gamma_0;p),
$$
\end{linenomath*}
for all $t\in \mathbb{R}$. The period of this solution is the minimum such $T$. It is convenient to parametrise points along a periodic solution by time $t$,
\begin{linenomath*}
$$
\gamma_t := X(t,\gamma_0;p)\in\mathbb{R}^2\quad\mbox{for any}\quad t \in [0,T),
$$
\end{linenomath*}
and define the corresponding cycle in the phase space in terms of $\gamma_t$,
\begin{linenomath*}
$$
\Gamma(p) = \{ \gamma_t: t\in[0,T) \} \subset \mathbb{R}^2.
$$
\end{linenomath*}

\begin{figure}[ht!]
\centering
\includegraphics[width=0.9\textwidth]{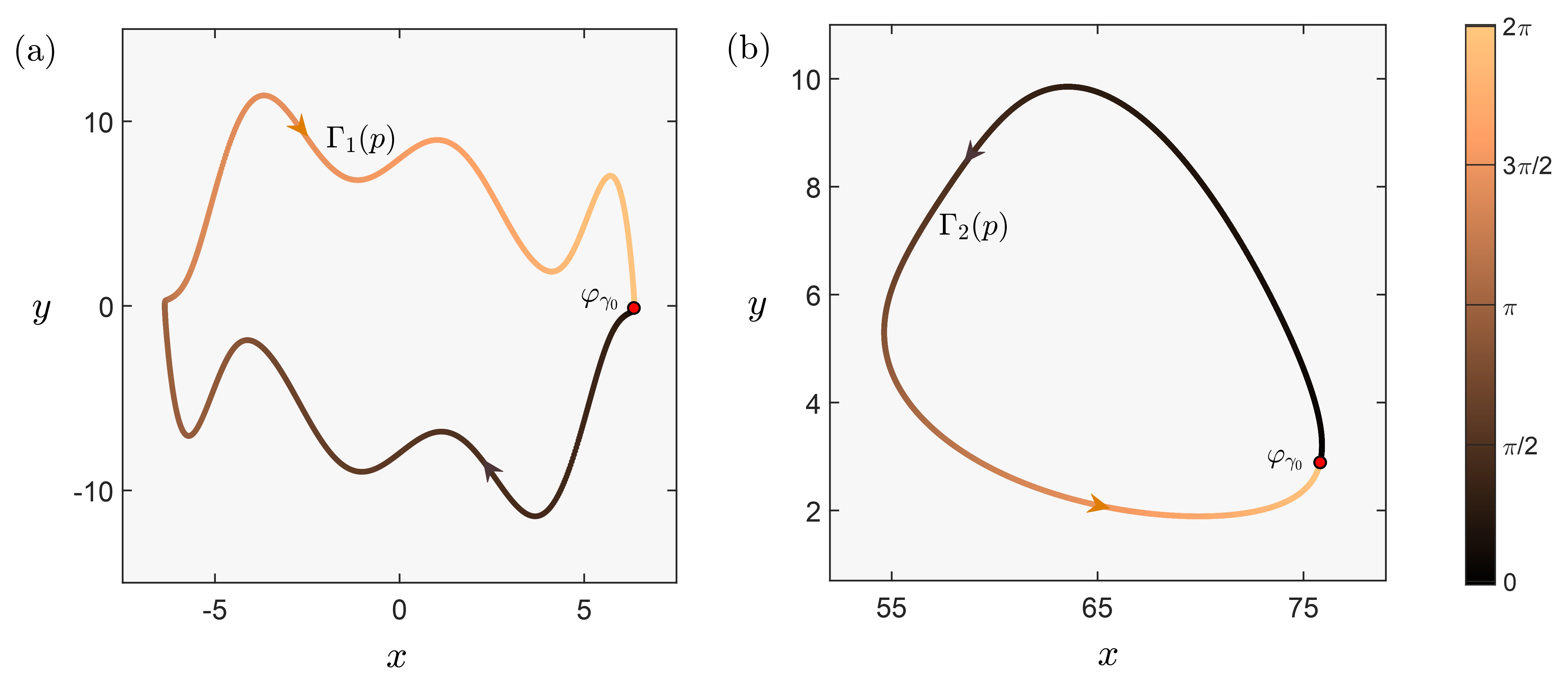}
\caption{Heat map distribution of the phases $\varphi_{\gamma_t} \in [0,2\pi)$ (defined in \eqref{phase_def}) corresponding to the points $\gamma_t$ on a limit cycle: For (a) $\Gamma_1(p)$ of the birhythmic van der Pol oscillator \eqref{vdp_with_xy}, where $p = \mu = 1.52$, and (b) $\Gamma_2(p)$ of the Decroly-Goldbeter glycolysis model \eqref{gly_model}, where $p = \sigma_i = 1.226$, starting at the initial phase $\varphi_{\gamma_0}$ (red filled circle).}
\label{fig: Figure_phase_of_LC}
\end{figure}

Following \citet{nakao2016phase}, we associate a {\em phase} $\varphi_{\gamma_t}$ with each point $\gamma_t$ along the cycle $\Gamma(p)$,  
\begin{linenomath*}
\begin{equation}
    \label{phase_def}
    \varphi_{\gamma_t} :=  \frac{2 \pi t}{T}\in [0, 2\pi).
\end{equation}
\end{linenomath*}
Here, we choose $\gamma_0$ with phase $ \varphi_{\gamma_0} = 0$ as follows: 
$\gamma_0 = (x_0, y_0) \in \Gamma_1(p)$, where $ x_0 = \max \{x: (x,y) \in \Gamma_1(p) \}$, for the birhythmic van der Pol oscillator, 
and similarly, $\gamma_0 = (x_0, y_0) \in \Gamma_2(p)$, where $ x_0 = \max \{x: (x,y) \in \Gamma_2(p) \}$, for the Decroly-Goldbeter glycolysis model. Fig.~\ref{fig: Figure_phase_of_LC} shows the phases for $\Gamma_1(p)$ in the van der Pol oscillator and $\Gamma_2(p)$ in the Decroly-Goldbeter glycolysis model.

\subsection{Partial basin instability of limit cycles \label{SSec: BI_limit_cycles}}

{\em Partial basin instability} for limit cycles~\citep{alkhayuon2021phase} is a simple geometric framework that uses global properties of the frozen autonomous system~\eqref{eq:auto_pbm_setting} to study RP-tipping. It compares the position of the base limit cycle at one point on a parameter path with the position of the basin of attraction of that limit cycle at other points on the parameter path. 

Consider the base limit cycle $\Gamma(p)$ and its basin of attraction $B(\Gamma,p)$. Also, consider a parameter path $\Delta_p$, on which $\Gamma(p)$ varies smoothly with $p$, meaning that $\Gamma(p)$ does not bifurcate along $\Delta_p$.
We say that $\Gamma(p)$ is {\em basin unstable} on $\Delta_p$ if there exist two points $p_1$ and $p_2$ in $\Delta_p$ such that $\Gamma (p_1)$ is not contained in the basin of attraction of $\Gamma (p_2)$,
\begin{linenomath*}
    \begin{equation}\label{eq:BI}
    \Gamma(p_1) \not\subset B(\Gamma,p_2). 
\end{equation}
\end{linenomath*}
Furthermore, $\Gamma(p)$ is {\em partially basin unstable} on $\Delta_p$ if, in addition to condition~\eqref{eq:BI}, for every two points $p_3, p_4 \in \Delta_p$, $\Gamma(p_3)$ has a nonempty intersection with $B(\Gamma,p_4)$. In other words, $\Gamma(p_3)$ never lies entirely outside of $B(\Gamma,p_4)$.
If there exist $\tilde{p}_3$ and $\tilde{p}_4$ such that  $\Gamma(\tilde{p}_3)$ has an empty intersection with $B(\Gamma,\tilde{p}_4)$, meaning that $\Gamma(\tilde{p}_3)$ lies entirely outside of $B(\Gamma,\tilde{p}_4)$, we say that $\Gamma(p)$ is {\em totally basin unstable} on $\Delta_p$.\footnote{Note that, there are indiscernible (or marginal) cases of basin instability that are not relevant to this work. For more details, see \cite[Sec. 4]{alkhayuon2021phase}.}

In the autonomous frozen system~\eqref{eq:auto_pbm_setting}, suppose that $\Gamma(p)$ is basin unstable on a parameter path $\Delta_p$. For any given $p_1$ and $p_2$ in $\Delta_p$, we will consider the set of all points  on the limit cycle $\Gamma(p_1)$ that lie outside ${B(\Gamma,p_2)}$,  
\begin{linenomath*}
\begin{equation*}
U(\Gamma,p_1,p_2) =
\{
\gamma_t: \gamma_t\, \in \Gamma(p_1)~\mbox{and}~\gamma_t\, \not\in B(\Gamma,p_2)
\}.
\end{equation*}
\end{linenomath*}
We will be interested in the corresponding set of phases of $\Gamma(p_1)$, namely
\begin{linenomath*}
\begin{equation*}
\Phi_U(\Gamma,p_1,p_2) =
\{
\varphi_{\gamma_t}: \gamma_t \in U(\Gamma,p_1,p_2)
\},
\end{equation*}
\end{linenomath*}
which we refer to as the $(p_1, p_2)$-basin unstable phases of $\Gamma(p)$, or basin unstable phases of $\Gamma(p)$ for short. 

Furthermore, for a given parameter path $\Delta_p$, we will consider the whole family of limit cycles $\Gamma(p)$ that vary smoothly with $p\in\Delta_p$,
\begin{linenomath*}
$$
G= \{\Gamma(p): p\in\Delta_p\}.
$$
\end{linenomath*}
For a fixed $p = p_1\in\Delta_p$, we define the {\em parameter region of basin instability} of $\Gamma(p_1)$ as follows:
\begin{linenomath*}
$$
BI(\Gamma,p_1) := \{ p_2\in\Delta_p: \Gamma(p_1)\not\subset B(\Gamma,p_2)\}.
$$
\end{linenomath*}
This is a set of all parameter values $p_2\in\Delta_p$ such that $\Gamma(p_1)$ is not contained in the basin of attraction of $\Gamma(p_2)$.

\subsection{Rate-induced phase-tipping from limit cycles\label{SSec: BI_implies_tipping}}

Consider the non-autonomous system~\eqref{eq:nonauto_pbm_setting} with sufficiently large $t_c$ so that it resembles the past limit system for small $t\ge 0$. We choose initial conditions $X_0$ on the base limit cycle $\Gamma(p^-)$  of the past limit system\footnote{Note that, $\Gamma(p^-)$ is not a limit cycle for system \eqref{eq:nonauto_pbm_setting} and $p(0)\neq p^-$, but close enough for large $t_c$, see Eqs.~\eqref{monotone_shift}~and~\eqref{nonmonotone_shift}.} and, for simplicity, we write 
$X(t,0,X_0) \ \to D$ as $t\to +\infty$
to denote $d[X(t,0,X_{0}), \, D] \to 0 ~~ \mbox{as} ~~ t\to +\infty$, where 
\begin{linenomath*}
    $$
 d[X,D] = \inf_{y\in{D}}\,\Vert X - y\Vert,
$$
\end{linenomath*}
is the Hausdorff semi-distance between a point $X$ and a compact set $D$.

We say that system~\eqref{eq:nonauto_pbm_setting} undergoes {\em irreversible} RP-tipping {from $\Gamma(p^-)$ for the time-dependent external inputs \eqref{monotone_shift} and \eqref{nonmonotone_shift}} if there exists a critical rate $r=r_c$ in the following sense. For all $0< r<r_c$ and every initial condition  $X_0 \in \Gamma(p^-)$, the solution $X(t,0,X_0)$ end-point tracks the base limit cycle, meaning that
\begin{linenomath*}
$$
X(t,0,X_0) \ \to \ \Gamma(p^+)~{\rm as}~t \ \to \ + \infty.
$$
\end{linenomath*}
However, for some $r>r_c$, there exist two points $X_{0,1}$ and $X_{0,2}$ on $\Gamma(p^-)$, such that the solution $X(t,0, X_{0,1})$ started from $X_{0,1}$ end-point tracks the base limit cycle, whereas $X(t, 0, X_{0,2})$ started from $X_{0,2}$ end-point tracks some alternative attractor $A(p)$ other than $\Gamma(p)$. To be precise,
\begin{linenomath*}
$$
X(t,0, X_{0,1})  \to  \Gamma(p^+)~{\rm as}~t  \to + \infty 
\quad\mbox{and}\quad
X(t,0, X_{0,2})  \to  A(p^+)~{\rm as}~t  \to + \infty.
$$
\end{linenomath*}
Now, to each initial condition $X_{0,2}\in\Gamma(p^-)$ that gives RP-tipping from $\Gamma(p^-)$, we assign a {\em (tipping) phase} which is the phase of that point as specified in Sec.~\ref{SSec: PhaseLC} and Fig.~\ref{fig: Figure_phase_of_LC}.

To make the connection between RP-tipping and partial basin instability, consider a partial basin unstable limit cycle $\Gamma(p)$ on a parameter path $\Delta_p$. 
Partial basin instability implies there are two points $p^-, p^+ \in \Delta_p$, such that
\begin{linenomath*}
\begin{equation}\label{eq:BI_implies_RP-tip}
   \Gamma(p^-) \not\subset B(\Gamma,p^+)
   \quad\mbox{and}\quad
   \Gamma(p^-) \cap B(\Gamma,p^+) \neq \emptyset\;. 
\end{equation}
\end{linenomath*}
Now, consider a monotone external input $p(rt)$, given by Eq.~\eqref{monotone_shift}, where $p(rt) \to p^\pm$ as $t\to \pm \infty$. One can show that the system can exhibit RP-tipping, as follows: 
\begin{itemize}
\item[1.] 
For a slow external input (vanishing $r$), \cite[Theroem III.1]{alkhayuon2018rate} shows that the system will exhibit tracking, i.e., there will be $\tilde{r}>0$ such that for every $X_{0} \in \Gamma(p^-)$ and $0 < r<\tilde{r}$, we have
\begin{linenomath*}
$$
X(t,0, X_{0}) \to \Gamma(p^+)~{\rm as}~t \to \infty.
$$
\end{linenomath*}
\item[2.]  
For very fast external input ($r\to\infty$), condition~\eqref{eq:BI_implies_RP-tip} implies that there exist two points $X_{0,1}, X_{0,2} \in \Gamma(p^-)$ such that the trajectory started from $X_{0,1}$ end-point tracks $\Gamma(p)$ and the trajectory started from $X_{0,2}$ does not,
\begin{linenomath*}
$$
X(t,0, X_{0,1})  \to  \Gamma(p^+)~{\rm as}~t \to + \infty
\quad\mbox{and}\quad X(t, 0, X_{0,2}) \not\to \Gamma(p^+) ~{\rm as}~t \to + \infty.
$$
\end{linenomath*}

\item [3.] By continuity, there must be at least one critical rate $r_c\geq\tilde{r}$ that gives  RP-tipping. 
\end{itemize}
In Section~\ref{Sec: Result}, we use the birhythmic van der Pol oscillator and Decroly-Goldbeter glycolysis model to illustrate RP-tipping.

\section{Results \label{Sec: Result}}
 
In this section, we analyse RP-tipping in both the birhythmic van der Pol and Decroly-Goldbeter glycolysis models. We start by identifying the regions of partial basin instability in the parameter space. Then, we illustrate RP-tipping in both models using tipping diagrams for the monotone and the non-monotone time-dependent external inputs. Unless stated otherwise, we will consider the stable limit cycle $\Gamma_1$ as the base state for the birhythmic van der Pol model, and $\Gamma_2$ as the base state for the Decroly-Goldbeter glycolysis model.

\subsection{Rate-induced Phase-tipping in the birhythmic van der Pol model \label{SSec: RP-tip_vanderPol}}

\subsubsection{Partial basin instability in the birhythmic van der Pol model \label{SSSec: Pbi in birhythmic model}}

\begin{figure}[ht!]
\centering
\includegraphics[width=0.9\textwidth]{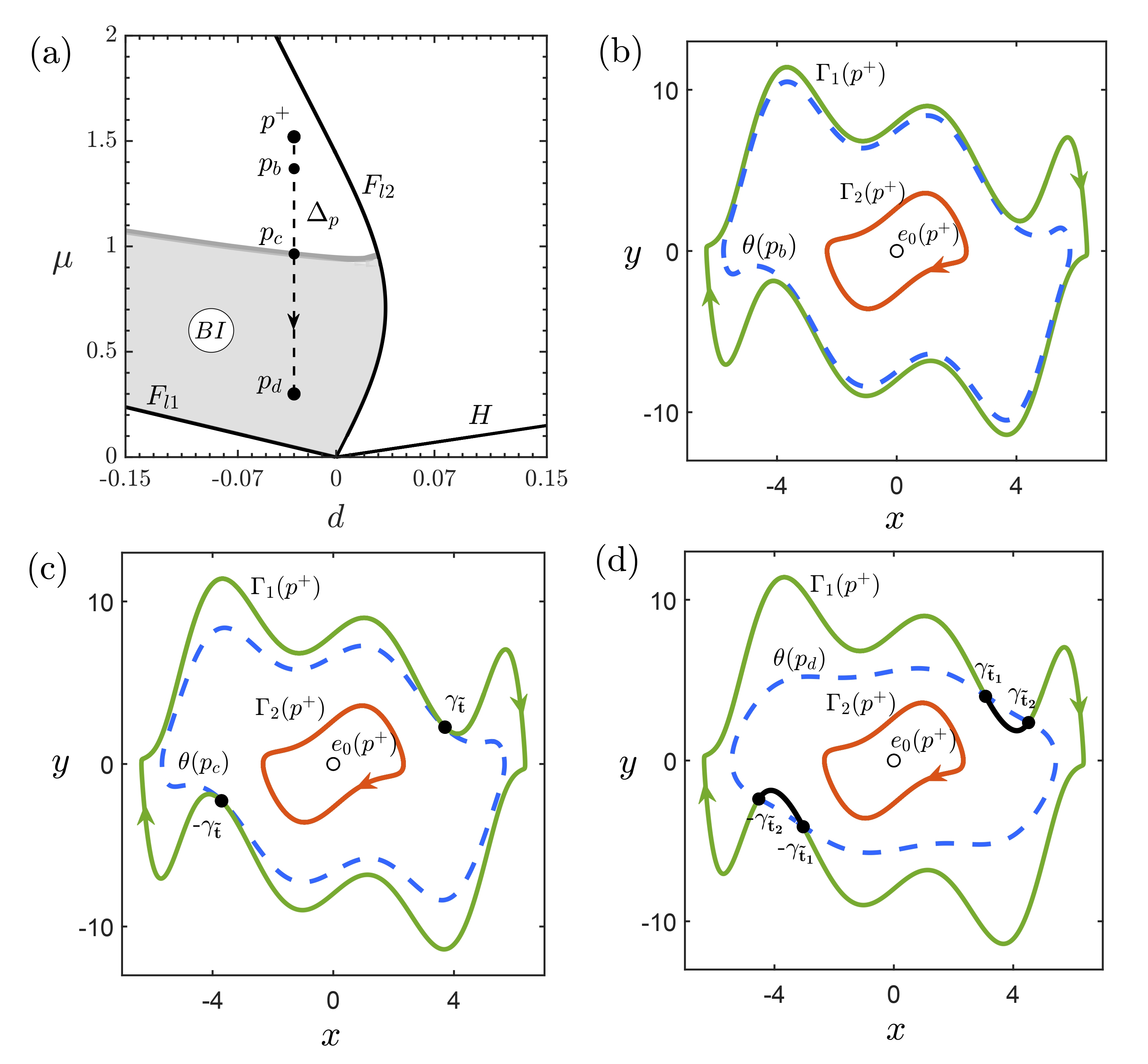}
\caption{(a) Two-parameter bifurcation diagram  ($d$ vs $\mu$) for the autonomous birhythmic van der Pol model (\ref{vdp_with_xy}), where the grey shaded area marks the region of partial basin instability $(BI)$ of the base limit cycle $\Gamma_1(p^+)$. The vertical dashed line between $p^+$ and $p_d$ represents a parameter path $\Delta_p$, which is used to simulate other sub-figures. (b)-(d) Phase portraits with two stable limit cycles $\Gamma_1$ (green) and $\Gamma_2$ (orange), and an unstable equilibrium point $e_0$ for a fixed $p^+ = (d,\mu^+)=(-0.03,1.52) \in \Delta_p$. In each phase portrait, the unstable limit cycle $\theta$ (blue dashed curve), which acts as a basin boundary between $\Gamma_1$ and $\Gamma_2$, is plotted for different parameter settings: (b) with no basin instability for $p_b = (-0.03,1.4) \in \Delta_p$, (c) with marginal basin instability at $\gamma_{\tilde{t}}$ (black dot) for $p_c = (-0.03,1) \in \Delta_p$ where $\tilde{t} \approx 6.11$, and (d) with partial basin instability between $\gamma_{\tilde{t}_1}$ $(\tilde{t}_1 \approx 5.9)$ and $\gamma_{\tilde{t}_2}$ $(\tilde{t}_2 \approx 6.51)$ for $p_d =(-0.03,0.3) \in \Delta_p$ (black curve on $\Gamma_1(p^+)$). The limit cycle $\Gamma_1(p^+)$ is also basin unstable at $-\gamma_{\tilde{t}}$ $(\tilde{t} \approx 2.61)$ in (c), as well as between $-\gamma_{\tilde{t}_1}$ $(\tilde{t}_1 \approx 2.4)$ and $-\gamma_{\tilde{t}_2}$ $(\tilde{t}_2 \approx 3.01)$ in (d), owing to the inherent symmetry of the model \eqref{vdp_with_xy}.}\label{fig: Figure_basin_ins_vdp}
\end{figure}

Consider the two-parameter bifurcation diagram Fig.~\ref{fig: Figure_bifn_vdp}(a) for the birhythmic van der Pol model. 
We are interested in the birhythmic region-I, where RP-tipping can occur from the base limit cycle $\Gamma_1$.  In Fig.~\ref{fig: Figure_basin_ins_vdp}(a), we augment the two-parameter bifurcation diagram Fig.~\ref{fig: Figure_bifn_vdp}(a) by adding the region of partial basin instability $BI$ for the base limit cycle $\Gamma_1$ at a particular parameter point $p^+$. 
In other words, the limit cycle $\Gamma_1(p^+)$ is partial basin unstable on any parameter path that crosses into the grey region $BI$. To illustrate this effect, we consider the parameter path $\Delta_p$, which does not cross any bifurcation curve, and consider three points along this path: $p_b$ in the white region,  $p_c$ on the boundary, and $p_d$ inside the grey region. Figs.~\ref{fig: Figure_basin_ins_vdp}(b)-\ref{fig: Figure_basin_ins_vdp}(d) show the stable limit cycles $\Gamma_{1,2}(p^+)$ in comparison with the unstable limit cycle $\theta(p)$, where $p = p_b$, $p_c,$ and $p_d$, respectively. 
Note that $\theta(p)$ represents the boundary of the basin of attraction of the base limit cycle.
One can see that in Fig.~\ref{fig: Figure_basin_ins_vdp}(b), the base limit cycle $\Gamma_1(p^+)$ does not intersect the basin boundary $\theta(p_b)$, indicating no basin instability.  However, in Fig.~\ref{fig: Figure_basin_ins_vdp}(d), there is a set of phases that are basin unstable (black) due to $\theta(p_d)$. The base limit cycle tangentially intersects the basin boundary $\theta(p_c)$ (see Fig.~\ref{fig: Figure_basin_ins_vdp}(c)) if the parameter path ends at point $p_c$; in such a case, $\Gamma_1(p^+)$ is marginally basin unstable~\citep{alkhayuon2021phase}.
From Section~\ref{SSec: BI_implies_tipping}, we expect RP-tipping from $\Gamma_1(p^+)$ if the external input varies on the parameter path $\Delta_p$. 

Furthermore, we point out that due to the symmetry of the model, the unstable limit cycle $\theta(p_c)$, in Fig.~\ref{fig: Figure_basin_ins_vdp}(c) tangentially intersects $\Gamma_1(p^+)$ at two points $\gamma_{\tilde{t}}$ and $-\gamma_{\tilde{t}}$, resulting in the possibility of basin instability at two locations and expanded further to give basin unstable phases between $\gamma_{\tilde{t}_1}$ ($-\gamma_{\tilde{t}_1}$) and $\gamma_{\tilde{t}_2}$ ($-\gamma_{\tilde{t}_2}$) in Fig.~\ref{fig: Figure_basin_ins_vdp}(d).

\begin{figure}[ht!]
\centering
\includegraphics[width=0.9\textwidth]{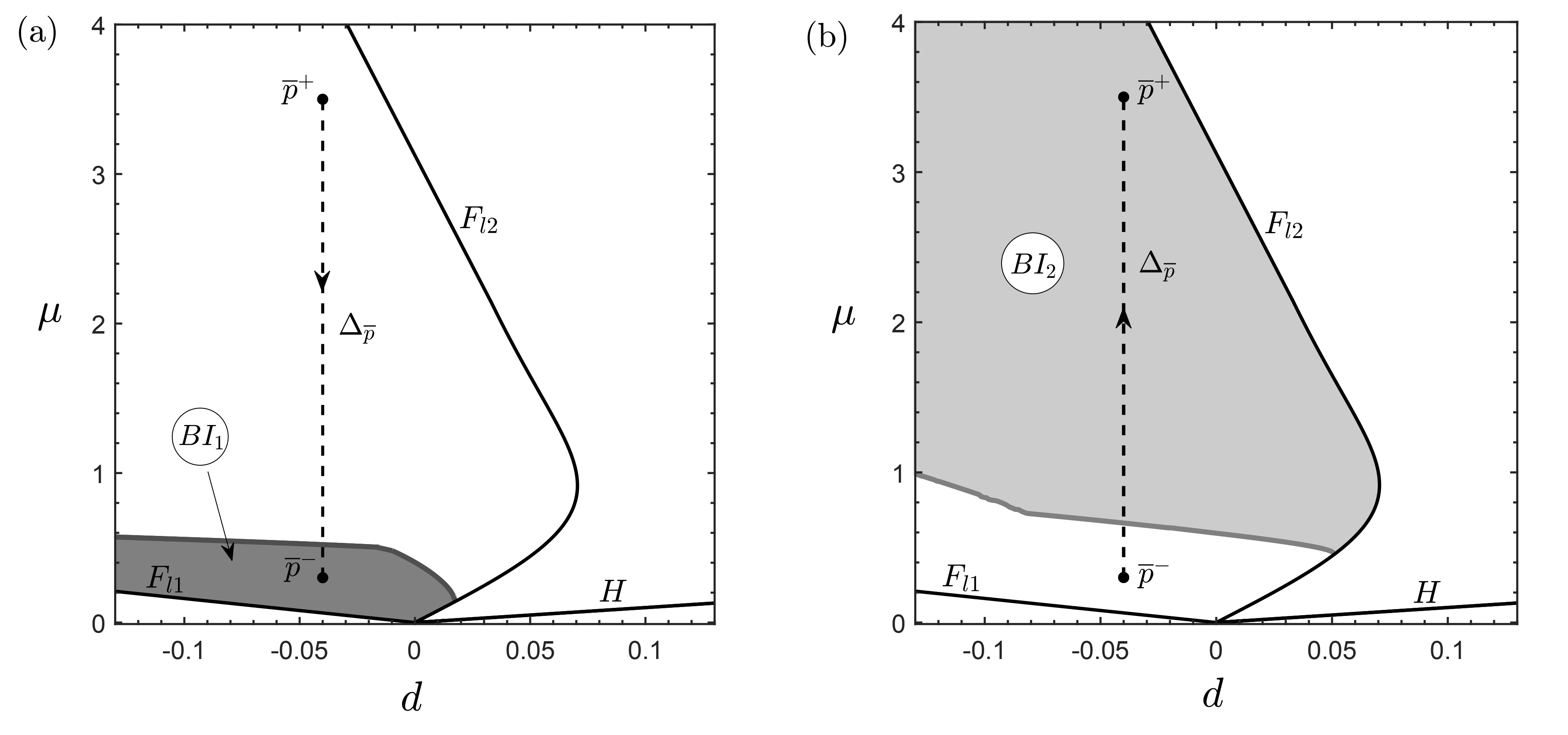}
 \caption{Two-parameter bifurcation diagrams of the birhythmic van der Pol model \eqref{vdp_with_xy} with the variation in the parameters $d$ vs. $\mu$. In (a), the dark grey shaded area denotes the region of partial basin instability $(BI_1)$ corresponding to the base state $\Gamma_2(\overline{p}^+)$ with the parameter path $\Delta_{\overline{p}}$ varied along the direction from $\overline{p}^+ = (-0.04, 3.5)$ to $\overline{p}^- = (-0.04, 0.3)$, whereas, the light grey region in (b) denotes the region of partial basin instability $(BI_2)$ correspond to the base state $\Gamma_1(\overline{p}^-)$ along with the parameter path $\Delta_{\overline{p}}$ varied from $\overline{p}^-$ to $\overline{p}^+$. Grey curves separating the grey and white regions in both figures represent the marginal basin unstable curves. The other parameter values are $\alpha = 0.0938$ and $\beta = 0.00194$.}
\label{fig: Figure_both_basin_ins}
\end{figure}

Considering the basin instability of $\Gamma_1(p)$ and the possibility of RP-tipping, one is tempted to look at the basin instability of the two stable limit cycles $\Gamma_1(p)$ and $\Gamma_2(p)$ for different points on the same parameter path. In Fig.~\ref{fig: Figure_both_basin_ins}, we consider two points $\overline{p}^-$ and $\overline{p}^+$ from a parameter path $\Delta_{\overline{p}}$ that is longer than previously considered $\Delta_p$. We then compute the basin instability region of $\Gamma_2(\overline{p}^+)$, denoted by $BI_1$ in Fig.~\ref{fig: Figure_both_basin_ins}(a), as well as, the basin instability region of $\Gamma_1(\overline{p}^-)$, denoted by $BI_2$ in Fig.~\ref{fig: Figure_both_basin_ins}(b). Although the $BI_1$ region is considerably small, having both basin instabilities on the same parameter path opens the possibility of alternating series of RP-tipping from $\Gamma_2(\overline{p})$ to $\Gamma_1(\overline{p})$ and vice versa if the external input varies non-monotonically on $\Delta_{\overline{p}}$. We illustrate this behaviour in Section~\ref{SSSec:  series of RP-tipping}.

\subsubsection{ Rate-induced Phase-tipping diagrams for the van der Pol model \label{SSSec: tipping dgms of vdP}}

\begin{figure}[ht!]
\centering
\includegraphics[width=0.9\textwidth]{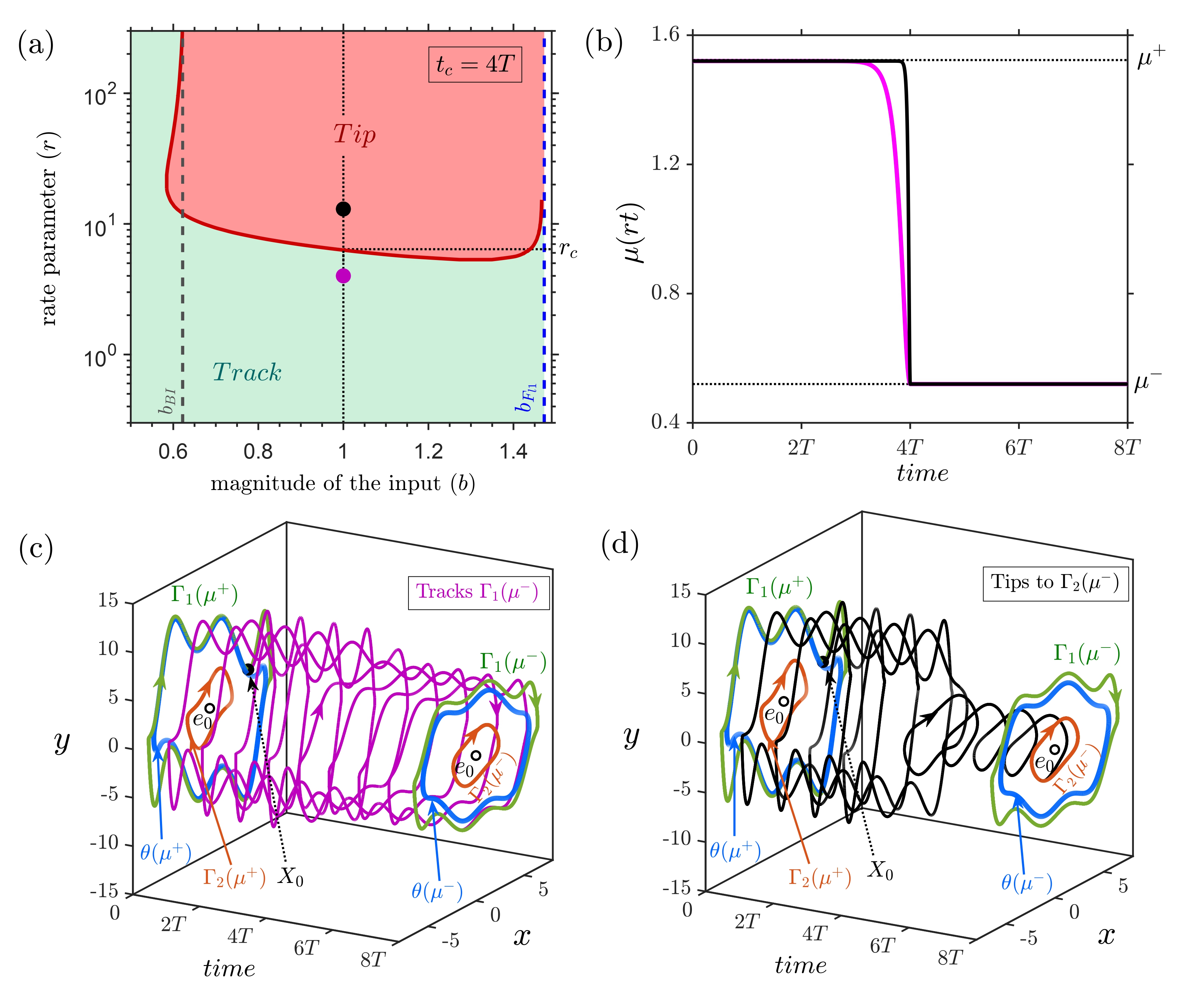}
\caption{(a) Half-bounded tipping diagram for the birhythmic van der Pol oscillator \eqref{vdp_with_xy} with variations in the parameters of monotone shift $\mu(rt)$ \eqref{monotone_shift}: magnitude ($b$) vs. rate ($r$) for the fixed value of $t_c = 4T$, where $T = 6.9979$ is period of the limit cycle $\Gamma_1(\mu^+)$. The red region (Tip) signifies the parameter space that results in RP-tipping, while the green region (Track) represents tracking. The red curve, which forms the boundary between the tipping and tracking regions, provides a critical rate $r_c$ of the monotone shift for a fixed magnitude $b$. The blue and grey vertical dashed lines indicate the fold point $b_{F_{l1}}$ and basin unstable phase $b_{BI}$, respectively. (b) The evolution of the monotone parameter shifts $\mu(rt)$  between $\mu^+$ and $\mu^-$ with respect to time for fixed $b$ and different values of the rate $r$ (same as the coordinate and colour scheme of the filled circles in (a)). (c)-(d) Phase portraits with trajectories originating from the initial point $X_0 = (4, 1.89)$ on $\Gamma_1(\mu^+)$ for the fixed magnitude $b = 1$ and different rates $r_1 = 4$ (magenta) and $r_2 = 13$ (black) of the monotone shift as pointed in (a). The trajectory tracks $\Gamma_1(\mu^-)$ for the rate $r_1$; whereas, for $r_2$ it tips to the limit cycle $\Gamma_2(\mu^-)$ (see the electronic supplementary material).} \label{fig: Figure_tipping_dgm_vdp_tanh}
\end{figure}

In this section, we use tipping diagrams to analyse RP-tipping in the birhythmic van der Pol model for both monotone and non-monotone external inputs. Throughout this section, the external input varies along the parameter path $\Delta_p$, presented in Fig.~\ref{fig: Figure_basin_ins_vdp}(a). In other words, we fix the parameter $d = -0.03$ and vary the input parameter $\mu$ following Eqs.~\eqref{monotone_shift} and \eqref{nonmonotone_shift} with $a = \mu^+$ and $b=\mu^+ - \mu^-$, we recall that $b$ is the magnitude of the external input. Rate-induced 
diagrams (or tipping diagrams for short) are two-parameter representations of tipping and tracking behaviour, i.e., the rate of change $r$ versus the magnitude of an external input $b$ \citep{o2020tipping}.

Fig.~\ref{fig: Figure_tipping_dgm_vdp_tanh}(a) presents 
a tipping diagram for 
model \eqref{vdp_with_xy}, 
with a monotone input $\mu(rt)$ \eqref{monotone_shift} that transitions from $\mu^+$ to $\mu^-$ along the path $\Delta_p$.
The figure is produced by solving the system for one initial condition $X_0 \approx (4, 1.89)$ on the base limit cycle $\Gamma_1(\mu^+)$ at time $t = 0$. The peak of the shift is fixed at time $t = t_c$, which is 4 times the period $(T)$ of the base limit cycle. 
The green region represents $r$ and $b$ values where the trajectory tracks the base limit cycle $\Gamma_1(\mu^-)$, whereas the red region represents $r$ and $b$ values where the trajectory tips to the alternative limit cycle $\Gamma_2(\mu^-)$. The red curve, which forms the boundary between these regions, indicates the critical rates $r_c$ for each value of $b$. 
Here, the tipping region is half-bounded. 
In the sense that, for each large enough magnitude $b$ there is a critical rate $r_c$ below which there is no tipping.
We point out that when $r$ is large, the boundary of the tipping region is asymptotic to $b_{BI} = 0.622$ (vertical grey dashed line), which corresponds to the basin unstable phase of an initial condition $X_0$.
On the other hand, we do not allow the external input to exceed the fold of limit cycles $b_{F_{l1}} = 1.4725$ (vertical blue dashed line).
To illustrate the tipping and tracking behaviour, we chose two points on the diagram that corresponds to the rates $r_1 = 4$ and $r_2 = 13$ (the magenta and black dots) for a fixed $b = 1$. Fig.~\ref{fig: Figure_tipping_dgm_vdp_tanh}(b) shows the external forcing in time concerning rates $r_1$ and $r_2$, and Fig.~\ref{fig: Figure_tipping_dgm_vdp_tanh}(c) shows the magenta trajectory with the rate $r_1$ is tracking $\Gamma_1(\mu^-)$,  while Fig.~\ref{fig: Figure_tipping_dgm_vdp_tanh}(d) shows the black trajectory with $r_2$ is tipping to the alternative limit cycle $\Gamma_2(\mu^-)$.

\begin{figure}[ht!]
\centering
\includegraphics[width=1\textwidth]{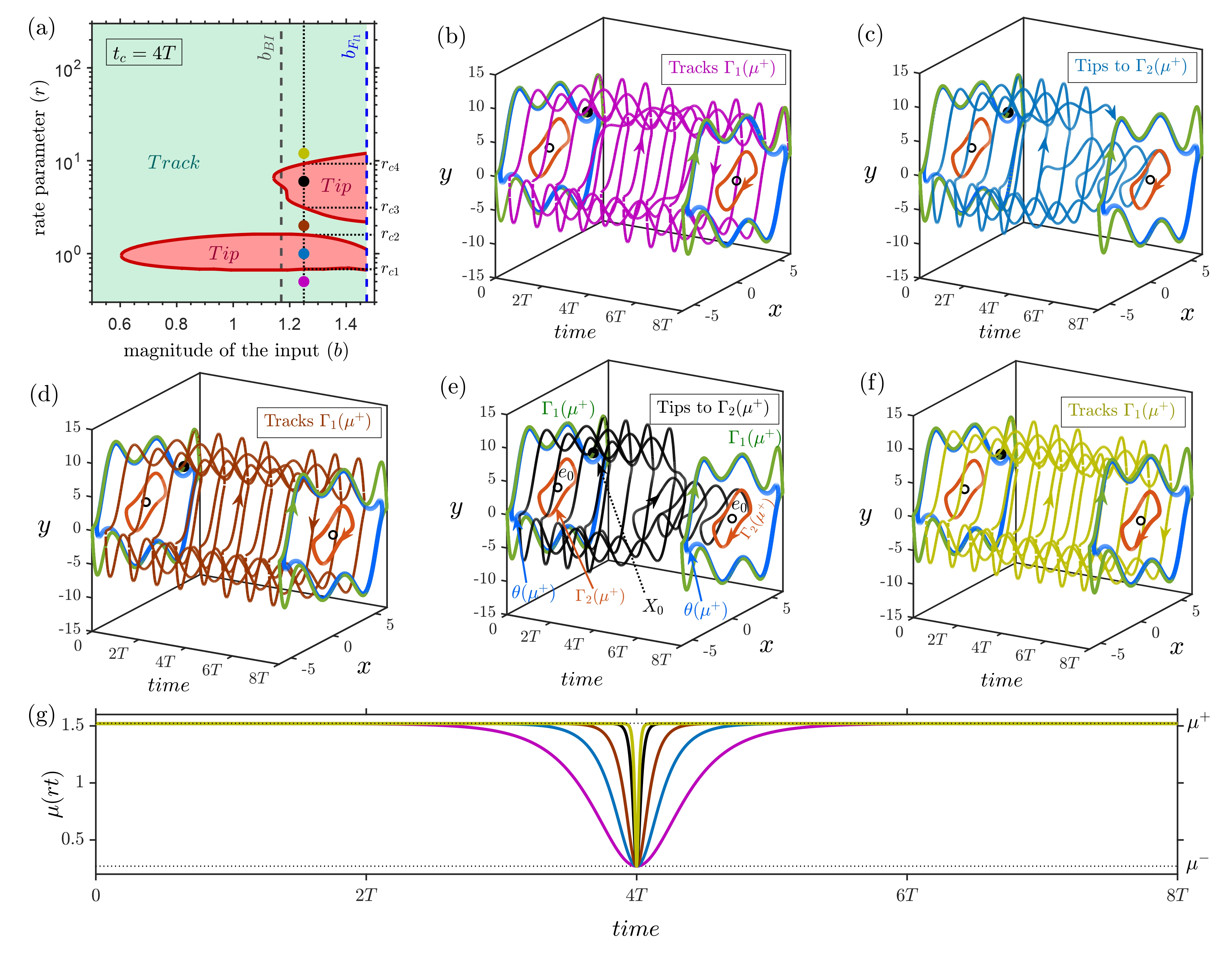}
\caption{(a) Bounded tipping diagram for the birhythmic van der Pol oscillator \eqref{vdp_with_xy} with variations in the parameters of non-monotone shift $\mu(rt)$ \eqref{nonmonotone_shift}: magnitude ($b$) vs. rate ($r$) for the fixed value of $t_c = 4T$, where $T$ is period of the limit cycle $\Gamma_1(\mu^+)$. The red region (Tip) signifies the parameter values that result in RP-tipping from $\Gamma_1(\mu^+)$ to $\Gamma_2(\mu^+)$, while the green region (Track) represents tracking to $\Gamma_1(\mu^+)$. The boundary curves (red) that separate the tipping and tracking regions provide critical rates $r_c$ of the non-monotone shift for a fixed magnitude $b$. The critical rates $r_{c1}, r_{c2}, r_{c3}$ and $r_{c4}$ correspond to the magnitude $b = 1.25$. The blue and grey vertical dashed lines indicate the fold point $b_{F_{l1}} = 1.4725$ and basin unstable phase $b_{BI} = 1.17$ corresponding to the initial condition $X_0$ in other sub-figures, respectively.  (b)-(f) Phase portraits with trajectories originating from the initial point $X_0 = (4.503, 2.33)$ on $\Gamma_1(\mu^+)$ for the fixed magnitude $b = 1.25$ and different rates $r_1 = 0.5$ (magenta), $r_2 = 1$ (blue), $r_3 = 2$ (brown), $r_4 = 6$ (black), and $r_5 = 12$ (yellow-green) as pointed in (a). For the rates $r_1$, $r_3$ and $r_5$, it is observed that trajectories track the limit cycle $\Gamma_1(\mu^+)$; whereas for the rates $r_2$ and $r_4$, trajectories tip to the alternative limit cycle $\Gamma_2(\mu^+)$ (see the electronic supplementary material). (g) The evolution of non-monotone parameter shifts $\mu(rt)$  with respect to time for fixed $b$ and different values of the rate $r$ (same as the coordinate and colour scheme of the filled circles in (a)).} 
 \label{fig: Figure_tipping_dgm_vdp_sech}
\end{figure}

Next, we consider model \eqref{vdp_with_xy} with a non-monotone external input from $\mu^+$ to $\mu^-$ and then back to $\mu^+$ following Eqn.~\eqref{nonmonotone_shift}, with the peak of the external input take place at time $t = t_c$, which is 4 times the period ($T$). 
Fig.~\ref{fig: Figure_tipping_dgm_vdp_sech}(a) presents a
tipping diagram with tipping regions bounded by tipping curves (red). 
In \citep{o2020tipping}, the authors refer to such bounded tipping regions as tipping tongues.
Similar to the monotone case, we do not allow the magnitude of the external input $b$ to exceed the fold of limit cycles $b_{F_{l1}} = 1.4725$ (blue dashed line). 

The tipping diagram delineates two distinct regions, red {\em(Tip)} and green {\em(Track)}, representing RP-tipping and tracking. The vertical grey dashed line indicates the basin unstable phase $b_{BI} = 1.17$ corresponds to the initial condition $X_0 \approx (4.503, 2.33)$ on the base limit cycle $\Gamma_1(\mu^+)$. The three-dimensional phase portraits in Figs.~\ref{fig: Figure_tipping_dgm_vdp_sech}(b)-\ref{fig: Figure_tipping_dgm_vdp_sech}(f) illustrate these tipping and tracking behaviours starting from $X_0$. For a fixed magnitude $b = 1.25$, five distinct rates are chosen (see Fig.~\ref{fig: Figure_tipping_dgm_vdp_sech} caption for rate values), each represented by a different colour in Fig.~\ref{fig: Figure_tipping_dgm_vdp_sech}(a). Among these, the rates $r_1$, $r_3$, and $r_5$ are located in tracking regions, while $r_2$ and $r_4$ are in tipping regions. As chosen rates in different regions, the trajectories in Figs.~\ref{fig: Figure_tipping_dgm_vdp_sech}(b), \ref{fig: Figure_tipping_dgm_vdp_sech}(d), and \ref{fig: Figure_tipping_dgm_vdp_sech}(f) track the same stable limit cycle $\Gamma_1(\mu^+)$, while the trajectories in Figs.~\ref{fig: Figure_tipping_dgm_vdp_sech}(c), and \ref{fig: Figure_tipping_dgm_vdp_sech}(e) tip to an alternative stable limit cycle $\Gamma_2(\mu^+)$. Fig.~\ref{fig: Figure_tipping_dgm_vdp_sech}(g) illustrates the non-monotone shifts corresponding to distinct rates. 

\begin{figure}[ht!]
\centering
\includegraphics[width=0.85\textwidth]{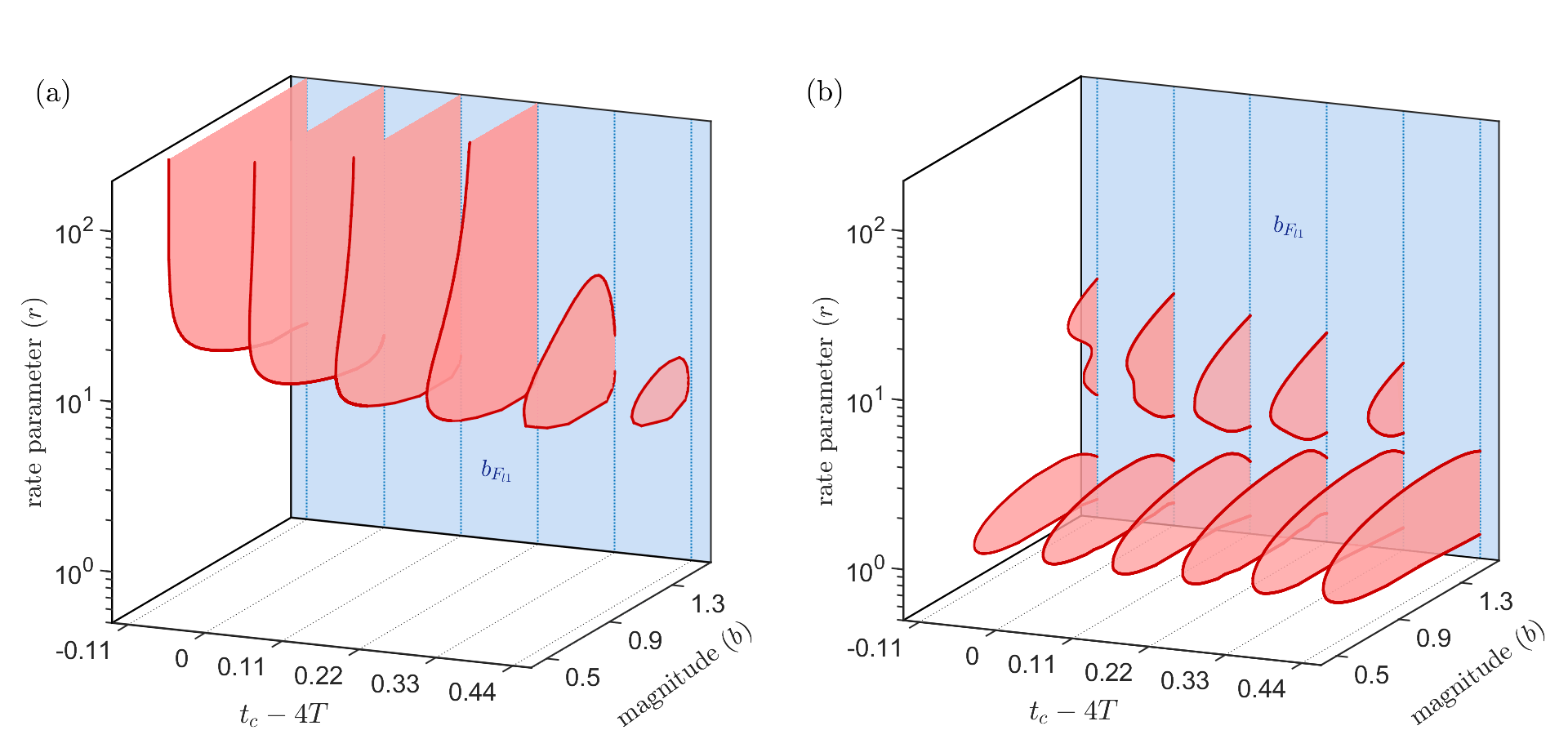}
\caption{Three-dimensional tipping diagrams of the birhythmic van der Pol model \eqref{vdp_with_xy} with variations in the parameters ($t_c$ vs. $b$ vs. $r$) for (a) monotone input \eqref{monotone_shift} and (b) non-monotone input \eqref{nonmonotone_shift}, for different values of $t_c$ starting from $4T - 0.11$ to $4T + 0.44$ with an increment of 0.11. The blue plane denotes the fold point $b_{F_{l1}}$.}
\label{fig: Figure_Multiple_tippg_curves_vdp}
\end{figure}

For a fixed value of $t_c = 4T$, which aligns with the phase of the initial condition $X_0$ of the trajectory, and for one trajectory, tipping diagrams reveal distinct characteristics of the tipping region: half-bounded for monotone input and fully bounded with two sub-regions for non-monotone input. This observation prompts an investigation into the relationship between the phases of the limit cycle and the tipping diagrams associated with both inputs described by \eqref{monotone_shift} and \eqref{nonmonotone_shift}. For a fixed initial condition $X_0$, one expects these tipping diagrams to vary while changing $t_c$. Also, here, we recall that our definition for the tipping phase is dependent on $t_c$. We find that changing the $t_c$ value alters the shape of rate-induced phase tipping diagrams, as illustrated in Fig.~\ref{fig: Figure_Multiple_tippg_curves_vdp}(a) for the monotone input. A gradual change in the value of $t_c$ by increments of 0.11 shows that the half-bounded tipping regions become fully bounded. Also, in the case of non-monotone input \eqref{nonmonotone_shift}, for variation in $t_c$ values, the impact of phase on the tipping diagrams can be seen in Fig.~\ref{fig: Figure_Multiple_tippg_curves_vdp}(b). It is observed that previously separated tipping regions merge, forming a single tipping region as $t_c$ varies.

\begin{figure}[ht!]
    \centering
    \includegraphics[width=0.9\textwidth]{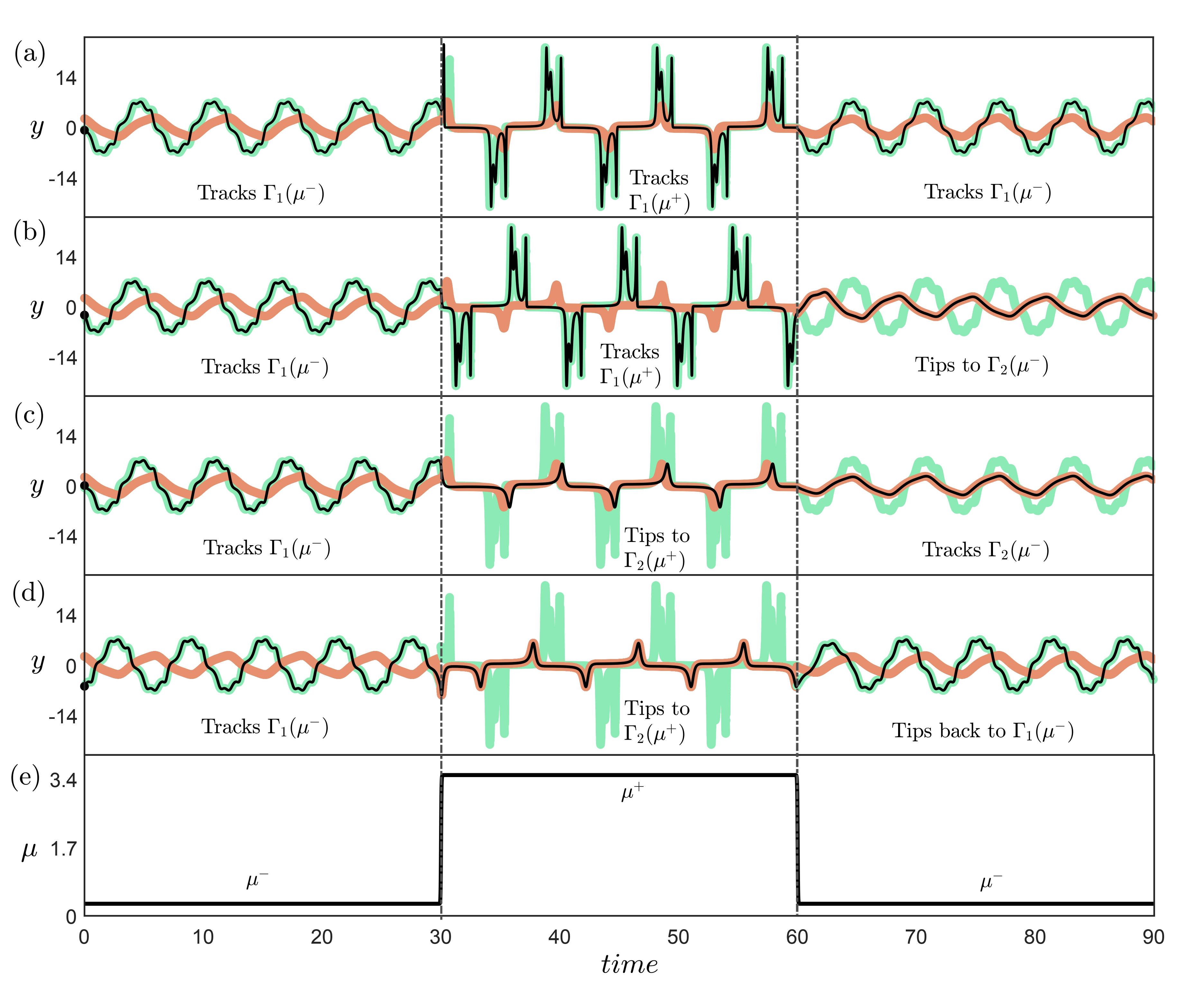}
    \caption{(a-d) Time series showing a series of tracking-tipping trajectories (black) of the birhythmic van der Pol model \eqref{vdp_with_xy} for the non-monotone input \eqref{impulse} at a single rate value $r = 27$ and different initial conditions on $\Gamma_1(\mu^-)$. These are shown along with the base stable limit cycles $\Gamma_1$ (green) and $\Gamma_2$ (orange) of the corresponding autonomous frozen system \eqref{vdp_with_xy}. In (a), the trajectory always tracks the base limit cycle $\Gamma_1$ with an initial condition $X_a = (6.45,-0.858)$. In (b), the trajectory tips to an alternative base state $\Gamma_2(\mu^-)$ at $t_{c2} = 60$ from $X_b = (5.77, -2.372)$. In (c), the trajectory from $X_c = (6.5726, -0.0013)$ tips to the alternative state $\Gamma_2(\mu^+)$ at $t_{c1} = 30$ and continues to track the same base state $\Gamma_2(\mu^-)$ despite of a change in the value of $\mu$. In (d), the trajectory from $X_d = (-2.7483, -5.9173)$ tips to $\Gamma_2(\mu^+)$ at $t_{c1} = 30$, tracks it until $t_{c2} = 60$, then tips back to $\Gamma_1(\mu^-)$. (e) The non-monotone time-dependent external input $\mu(rt)$ \eqref{impulse} in the form of an impulse.}
    \label{fig: series_of_tipping}
\end{figure}

\subsubsection{Series of rate-induced phase-tipping between limit cycles \label{SSSec: series of RP-tipping}}

In Section~\ref{SSSec: Pbi in birhythmic model}, we pointed out that there is a possibility of a series of alternating tipping between $\Gamma_1(p)$ and $\Gamma_2(p)$ in the van der Pol model. 
In this section, we illustrate this behaviour using a new non-monotone external input in the form of impulse  \citep{hasan2023rate}:
\begin{linenomath*}
\begin{equation}
\label{impulse}
    \mu(rt) = \mu^- + b \left [ \dfrac{\tanh(r(t-t_{c1})) - \tanh(r(t-t_{c2}))}{2} \right].
\end{equation}
\end{linenomath*}
In this non-monotone input, $\mu(rt) \approx \mu^+= \mu^- + b$ for a time interval of length $t_{c2} - t_{c1}$, which allows the system to settle into the new regime before shifting it back to $\mu(rt) \approx \mu^-$, see Fig.~\ref{fig: series_of_tipping}(e). In Fig.~\ref{fig: series_of_tipping}, we vary the input parameter along the parameter path $\Delta_{\overline{p}}$ shown in Fig.~\ref{fig: Figure_both_basin_ins}, where we fix the parameter $d = -0.04$ and vary the parameter $\mu(rt)$ between $\mu^- = 0.3$ and $\mu^+ = 3.5$ according to \eqref{impulse}. We also fix the rate parameter $r = 27$ and illustrate possible tipping and tracking combinations. The non-autonomous trajectories are depicted in solid black; green and orange trajectories represent base limit cycles $\Gamma_1$ and $\Gamma_2$, respectively, for fixed in-time input $\mu$. The 4 trajectories in Figs.~\ref{fig: series_of_tipping}(a)-\ref{fig: series_of_tipping}(d) start at different initial conditions $X_a,X_b,X_c$ and $X_d$ all on the limit cycle $\Gamma_1(\mu^-)$, at time $t = 0$. In Fig.~\ref{fig: series_of_tipping}(a), the trajectory tracks $\Gamma_1$ all the way to the end. In Fig.~\ref{fig: series_of_tipping}(b), the trajectory tracks initially but then tips to $\Gamma_2(\mu^-)$ after the second switch. In Fig.~\ref{fig: series_of_tipping}(c), the trajectory tips at the first switch to $\Gamma_2(\mu^+)$ and keeps tracking $\Gamma_2$ afterwords. Finally, in Fig.~\ref{fig: series_of_tipping}(d), the trajectory tips twice, from $\Gamma_1(\mu^-)$ to $\Gamma_2(\mu^+)$ at the first switch and back to $\Gamma_1(\mu^-)$ at the second switch.

\subsection{Rate-induced phase-tipping in the Decroly-Goldbeter glycolysis model \label{SSec: RP-tip_DGg}}

\subsubsection{Partial basin instability in the Decroly-Goldbeter glycolysis model \label{SSSec: Pbi in glycolysis model}}
\begin{figure}[ht!]
\centering
\includegraphics[width=0.9\textwidth]{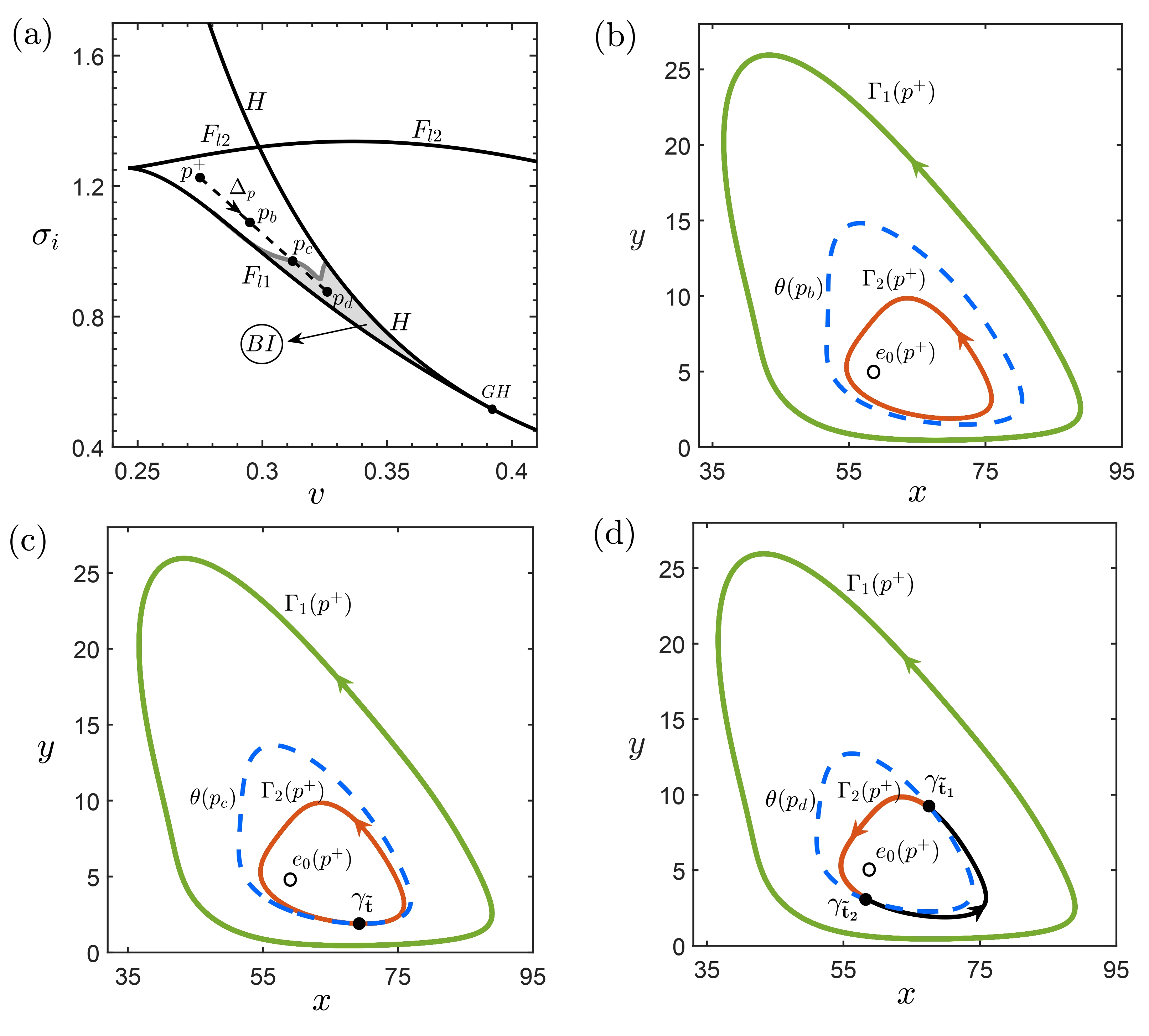}
\caption{(a) Two-parameter bifurcation diagram ($v$ vs $\sigma_i$) for the autonomous glycolysis model \eqref{gly_model}, where the grey shaded area marks the region of partial basin instability $(BI)$ corresponding to the base limit cycle $\Gamma_2(p^+)$. (b)-(d) Phase portraits with two stable limit cycles $\Gamma_1$ (green) and $\Gamma_2$ (orange), and an unstable equilibrium point $e_0$ for a fixed $p^+ = (v^+,\sigma_i^+)=(0.275,1.226)$. In each phase portrait, the unstable limit cycle $\theta$ (blue dashed curve), which acts as a basin boundary between $\Gamma_1$ and $\Gamma_2$, is plotted for different parameter settings along the parameter path $\Delta_p$: (b) with no basin instability for $p_b = (0.295,1.0887)$, (c) with marginal basin instability at $\gamma_{\tilde{t}}$ (black dot) for $p_c = (0.312,0.9702)$ where $\tilde{t} \approx 249.8$, and (d) with partial basin instability between $\gamma_{\tilde{t}_1}$ $(\tilde{t}_1 \approx 32.17)$ and $\gamma_{\tilde{t}_2}$ $(\tilde{t}_2 \approx 179.31)$ for $p_d =(0.326,0.876)$ (black curve on $\Gamma_2(p^+)$).}
\label{fig: Figure_basin_ins_gly}
\end{figure}

Here, we study partial basin instability in the model \eqref{gly_model} of birhythmic glycolytic oscillations. Consider the two-parameter bifurcation diagram in Fig.~\ref{fig: Figure_bifn_gly}(a). We are interested in the birhythmic region-I, where RP-tipping can occur from the base limit cycle. In Fig.~\ref{fig: Figure_basin_ins_gly}(a), we augment the two-parameter bifurcation diagram by adding the region of partial basin instability $BI$ for the base limit cycle at a particular parameter point $p^+$. In contrast to the van der Pol model \eqref{vdp_with_xy}, in this model, the limit cycle $\Gamma_2(p^+)$ is basin unstable on any parameter path that crosses into the grey region $BI$. To illustrate this effect, consider the parameter path $\Delta_p$, which does not cross any bifurcation curve, and consider three points along this path: $p_b$ in the white region, $p_c$ on the boundary, and $p_d$ in the grey region. This path is strategically chosen along the diagonal connecting points $p^+$ and $p_d$ to induce basin instability. The corresponding Figs.~\ref{fig: Figure_basin_ins_gly}(b),~\ref{fig: Figure_basin_ins_gly}(c), and \ref{fig: Figure_basin_ins_gly}(d) show the stable limit cycles $\Gamma_{1,2}(p^+)$ in comparison with the unstable limit cycle $\theta(p)$, where $p = p_b$, $p_c,$ and $p_d$. $\theta(p)$ represents the boundary of the basin of attraction of the base limit cycle. One can see that in Fig.~\ref{fig: Figure_basin_ins_gly}(b), the base limit cycle does not intersect the basin boundary and, hence, has no basin instability. However, in Fig.~\ref{fig: Figure_basin_ins_gly}(d), there is a set of basin unstable phases (black). The base limit cycle tangentially intersects the basin boundary if the parameter path ends at point $p_c$ (see Fig.~\ref{fig: Figure_basin_ins_gly}(c)); in such a case, the base limit cycle is marginally basin unstable~\citep{alkhayuon2021phase}. Again, from Section~\ref{SSec: BI_implies_tipping}, we expect RP-tipping from $\Gamma_2(p^+)$ if the external input varies on the parameter path $\Delta_p$.

\subsubsection{Rate-induced Phase-tipping diagrams of the Decroly-Goldbeter glycolysis model \label{SSSec: tipping dgms of gly}}


\begin{figure}[ht!]
\centering
\includegraphics[width=0.9\textwidth]{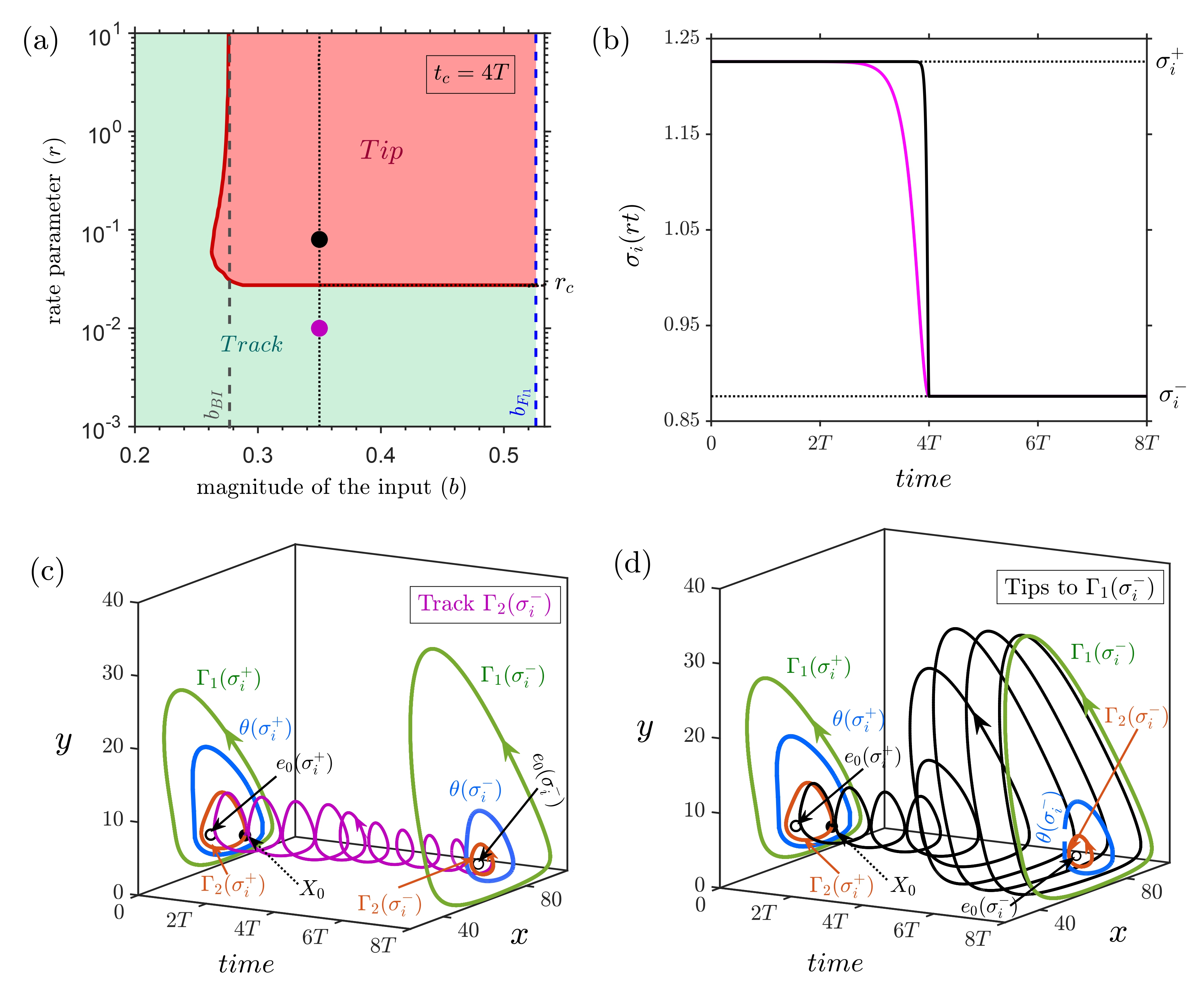}
\caption{(a) Half-bounded tipping diagram for the Decroly-Goldbeter glycolysis model \eqref{gly_model} with variations in the parameters of monotone shift $\sigma_{i}(rt)$ \eqref{monotone_shift}: magnitude ($b$) vs. rate ($r$) for the fixed value of $t_c = 4T$, where $T=308.266$ is the period of $\Gamma_2(\sigma^+_i)$. The red region (Tip) signifies the parameter values that result in RP-tipping, while the green region (Track) represents tracking. The horizontal black dotted line corresponds to the critical rate $r_c$ on the boundary curve for a fixed magnitude $b$. The blue and grey dashed lines indicate the fold point $b_{F_{l1}} = 0.526$ and basin unstable phase $b_{BI} = 0.2762$, respectively. (b) The evolution of monotone parameter shifts $\sigma_{i}(rt)$ between $\sigma_i^+$ and $\sigma_i^-$ with respect to time, for fixed $b$ and different values of the rate $r$ (same as the coordinate and colour scheme of the filled circles in (a)). (c)-(e) Phase portraits with trajectories originating from the initial point $X_0 = (75.71, 2.76)$ on $\Gamma_2(\sigma_i^+)$, for the fixed magnitude $b = 0.35$ and different rates $r_1 = 0.01$ (magenta) and $r_2 = 0.08$ (black) of the monotone shift as pointed in (a). The trajectory tracks the limit cycle $\Gamma_2(\sigma_i^+)$ for the rate $r_1$; whereas for the rate $r_2$, the trajectory tips to the limit cycle $\Gamma_1(\sigma_i^-)$ (see the electronic supplementary material).}
\label{fig: Figure_tipping_dgm_gly_tanh}
\end{figure}

Similar to the van der Pol model, we investigate RP-tipping in the model~\eqref{gly_model} using the tipping diagrams for both monotone and non-monotone external inputs. The external input varies along the parameter path $\Delta_p$, as shown in Fig.~\ref{fig: Figure_basin_ins_gly}(a). Specifically, we vary the parameter $\sigma_i$ between $\sigma_i^+$ and $\sigma_i^-$ along the diagonal parameter path $v \approx (-\sigma_i+3.11)/6.86$ following Eqs.~\eqref{monotone_shift} and \eqref{nonmonotone_shift}, with $a = \sigma_i^+$ and $b = \sigma_i^+ - \sigma_i^-$, where $b$ represents the magnitude of the external input.

First, we examine the model \eqref{gly_model} with a monotone shift \eqref{monotone_shift} of the external input $\sigma_i(rt)$ where the parameter transitions from $\sigma_i^+$ to $\sigma_i^-$. Fig.~\ref{fig: Figure_tipping_dgm_gly_tanh}(a) displays a tipping region for a fixed value of $t_c$, which is half-bounded \citep{o2020tipping}. Specifically, as the rate $r$ increases, the boundary of the tipping region asymptotically approaches $b_{BI} = 0.2762$ (grey dashed line), corresponding to the basin instability of the initial condition $X_0 \approx (75.706, 2.76)$. The peak of the shift is set at $t = t_c$, which is four times the period $(T)$ of the limit cycle. The tipping diagram includes a fold point $b_{F_{l1}} = 0.526$ (blue dashed line) beyond which the system exhibits B-tipping. The green region indicates values of $r$ and $b$ where the trajectory follows the base limit cycle $\Gamma_2(\sigma_i^-)$, whereas the red region indicates values of $r$ and $b$ where the trajectory tips to an alternative limit cycle $\Gamma_1(\sigma_i^-)$. The boundary curve (red) delineating these regions denotes the critical rates $r_c$ for a fixed value of $b$. To illustrate the tipping and tracking behaviour, two points on the diagram for rates $r_1$ and $r_2$ (the magenta and black dots) were selected for the fixed $b = 0.35$. Fig.~\ref{fig: Figure_tipping_dgm_gly_tanh}(b) shows the external forcing corresponding to both rates $r_1$ and $r_2$. Figs.~\ref{fig: Figure_tipping_dgm_gly_tanh}(c) and \ref{fig: Figure_tipping_dgm_gly_tanh}(d) demonstrate that the magenta trajectory with rate $r_1$ tracks the base limit cycle $\Gamma_2(\sigma_i^-)$, while the black trajectory with rate $r_2$ tips to the alternative limit cycle $\Gamma_1(\sigma_i^-)$, respectively.

\begin{figure}[ht!]
\centering
\includegraphics[width=0.95\textwidth]{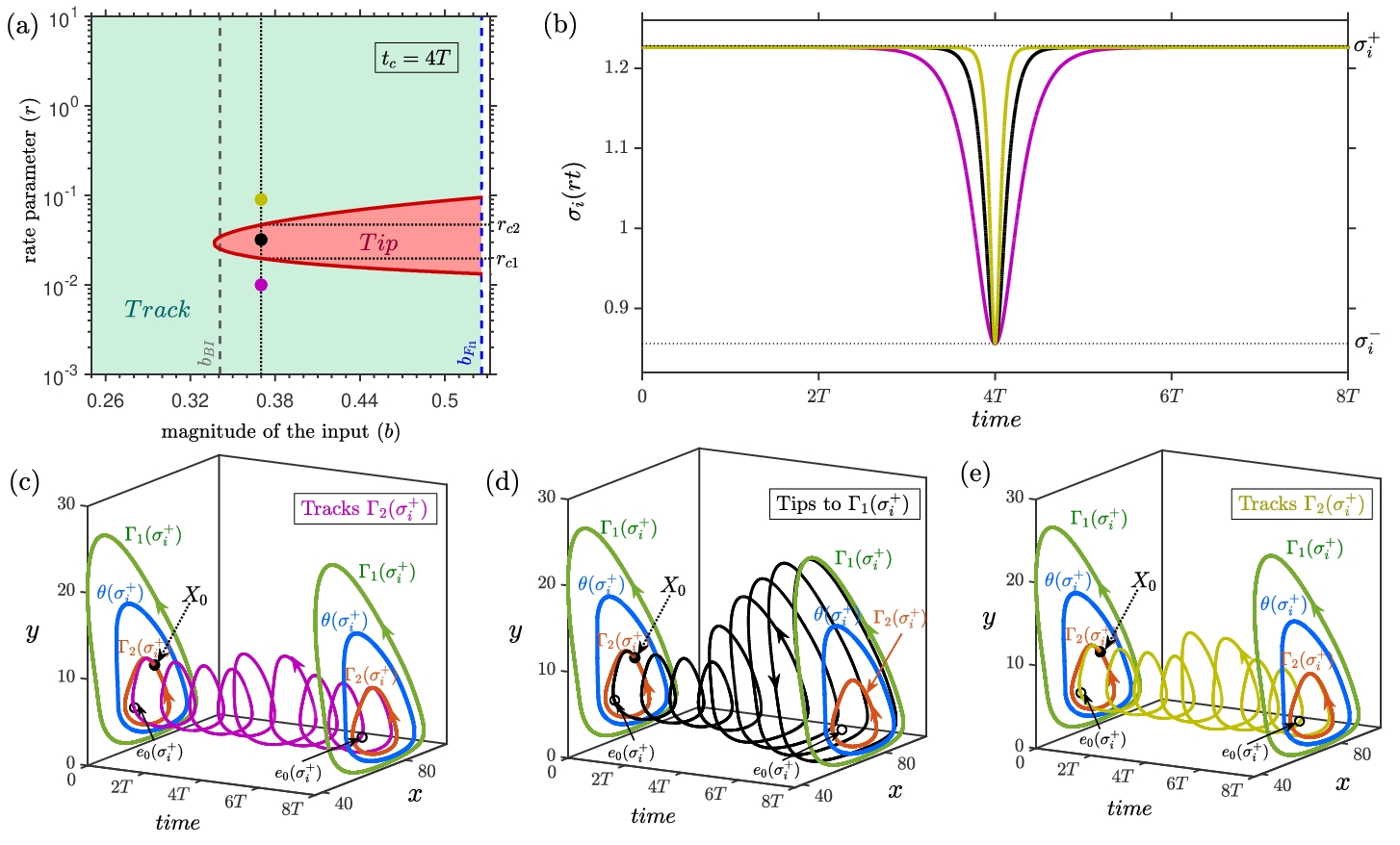}
\caption{(a) Bounded tipping diagram for the Decroly-Goldbeter glycolysis model \eqref{gly_model} with variations in the parameters of non-monotone shift $\sigma_{i}(rt)$ \eqref{nonmonotone_shift}: magnitude ($b$) vs. rate ($r$) for the fixed value of $t_c = 4T$, where $T$ is the period of the base limit cycle $\Gamma_2(\sigma^+_i)$. The red and green regions denote RP-tipping and tracking, separated by a boundary curve (red). $r_{c1}$ and $r_{c2}$ represent critical rates correspond to the magnitude $b = 0.37$. The blue and grey dashed lines indicate the fold point $b_{F_{l1}} = 0.526$ and basin unstable phase $b_{BI} = 0.341$, respectively. (b) The evolution of non-monotone parameter shifts $\sigma_{i}(rt)$ with respect to time for fixed $b$ and different values of the rate $r$ (same as the coordinate and colour scheme of the filled circles in (a)). (c)-(e) Phase portraits with trajectories originating from the initial point $X_0 = (68.91, 8.61)$ on $\Gamma_2(\sigma_i^+)$, for the fixed magnitude $b = 0.37$ and different rates $r_1 = 0.01$ (magenta), $r_2 = 0.032$ (black), and $r_3 = 0.09$ (grey) of the non-monotone shift as pointed in (a). For the rates $r_1$ and $r_3$, it is observed that the trajectory tracks the limit cycle $\Gamma_2(\sigma_i^+)$; whereas for the rate $r_2$, the trajectory tips to the limit cycle $\Gamma_1(\sigma_i^+)$ (see the electronic supplementary material).} \label{fig: Figure_tipping_dgm_gly_sech}
\end{figure}

Next, we examine model \eqref{gly_model} with a non-monotone shift \eqref{nonmonotone_shift} of the external input $\sigma_i(rt)$ where the parameter transitions from $\sigma_i^+$ to $\sigma_i^-$ and then returns to $\sigma_i^+$. Fig.~\ref{fig: Figure_tipping_dgm_gly_sech}(a) displays a tipping tongue for a fixed value of $t_c$. This figure is generated by solving the system from an initial condition $X_0 \approx (68.91, 8.614)$ on the base limit cycle $\Gamma_2(\sigma_i^+)$ at time $t=0$. The basin unstable point $b_{BI} = 0.341$ corresponds to $X_0$ and is marked with a vertical grey dashed line. The fold point $b_{F_{l1}} = 0.526$ is indicated by a blue dashed line, beyond which B-tipping occurs. The green region indicates values of $r$ and $b$ where the trajectory follows the base limit cycle $\Gamma_2(\sigma_i^+)$, whereas the red region indicates values of $r$ and $b$ where the trajectory tips to an alternative limit cycle $\Gamma_1(\sigma_i^+)$. Figs.~\ref{fig: Figure_tipping_dgm_gly_sech}(c)-\ref{fig: Figure_tipping_dgm_gly_sech}(e) illustrate the tipping and tracking behaviour in three-dimension. Three distinct rates are selected for a fixed magnitude $b = 0.37$, each shown in a different colour as mentioned in Fig.~\ref{fig: Figure_tipping_dgm_gly_sech}(a). The rates $r_1$ and $r_3$ lie within tracking regions, while $r_2$ lies in the tipping region. Consequently, the trajectories in Fig.~\ref{fig: Figure_tipping_dgm_gly_sech}(c) and \ref{fig: Figure_tipping_dgm_gly_sech}(f) follow the same stable limit cycle $\Gamma_2(\sigma_i^+)$, while the trajectory in Fig.~\ref{fig: Figure_tipping_dgm_gly_sech}(d) tips to an alternative stable limit cycle $\Gamma_1(\sigma_i^+)$. Fig.~\ref{fig: Figure_tipping_dgm_gly_sech}(b) shows non-monotone shifts corresponding to different rates. 

\begin{figure}[ht!]
    \centering
    \includegraphics[width=0.85\textwidth]{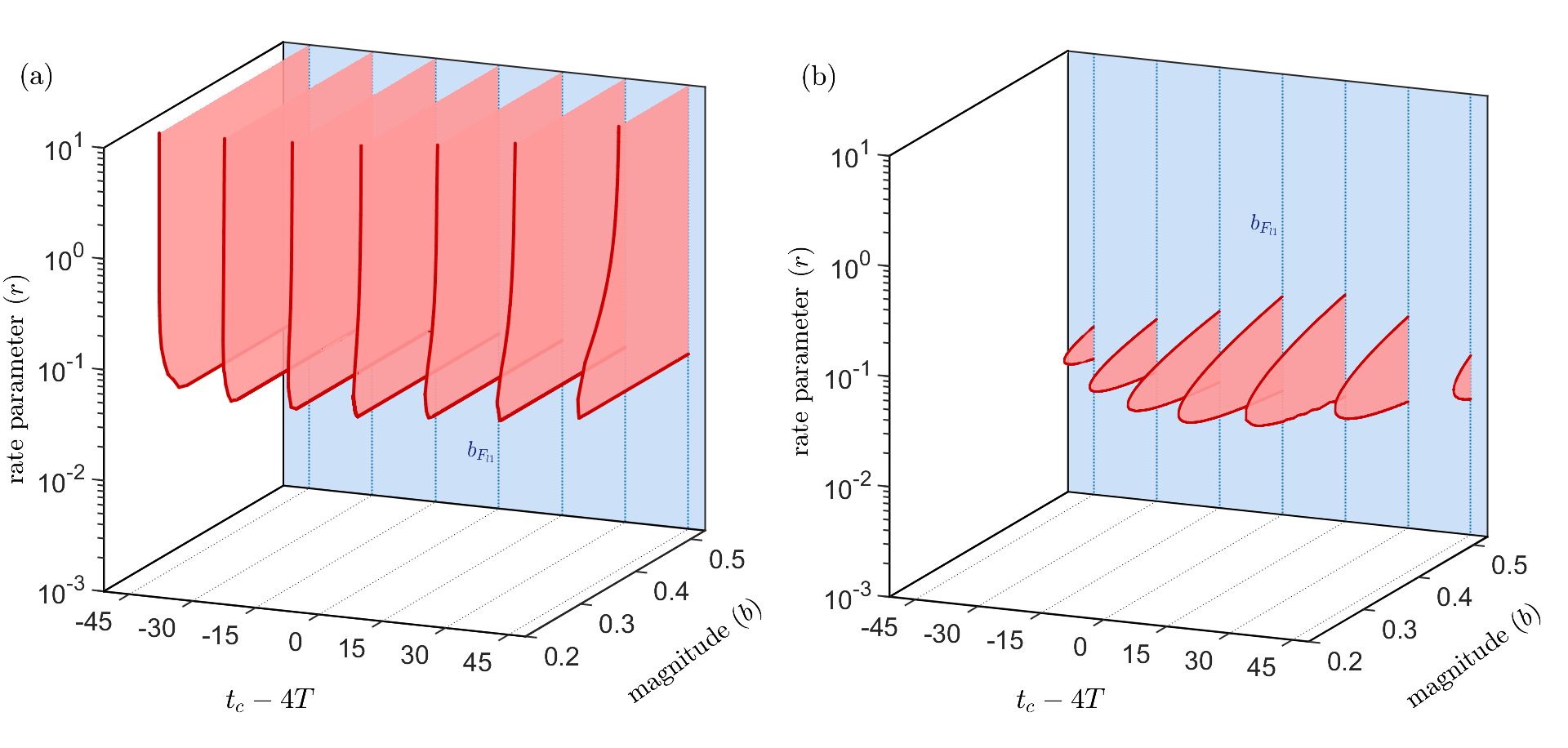}
    \caption{The set of three-dimensional tipping diagrams of the Decroly-Goldbeter glycolysis model \eqref{gly_model} with variations in the parameters ($t_c$ vs. $b$ vs. $r$) of (a) monotone input \eqref{monotone_shift} and (b) non-monotone input \eqref{nonmonotone_shift} for different values of $t_c$ starting from $4T - 45$ to $4T + 45$ in the interval of 15. The blue plane denotes the fold point $b_{F_{l1}}$.}
    \label{fig: Figure_Multiple_tippg_curves_gly}
\end{figure}


The impact of varying $t_c$ values on the tipping diagrams is shown in Fig.~\ref{fig: Figure_Multiple_tippg_curves_gly}(a) and Fig.~\ref{fig: Figure_Multiple_tippg_curves_gly}(b) with the incremental value of 15. The half-bounded tipping regions correspond to monotone input \eqref{monotone_shift} and gradually decrease in size as $t_c$ varies. In the case of non-monotone input, unlike the tipping diagram correspond to the van der Pol model in Fig.~\ref{fig: Figure_tipping_dgm_vdp_sech}(a), for a fixed value of $t_c = 4T$, the tipping diagram in glycolysis model with non-monotone shift \eqref{nonmonotone_shift} has a single bounded region. Upon varying the values of $t_c$, Fig.~\ref{fig: Figure_Multiple_tippg_curves_gly}(b) illustrates the change in the size of the tipping region.

\subsection{General relationship between phase and pace (rate) \label{SSec: pace vs phase}}

\begin{figure}[ht!]
    \centering
    \includegraphics[width=0.9\textwidth]{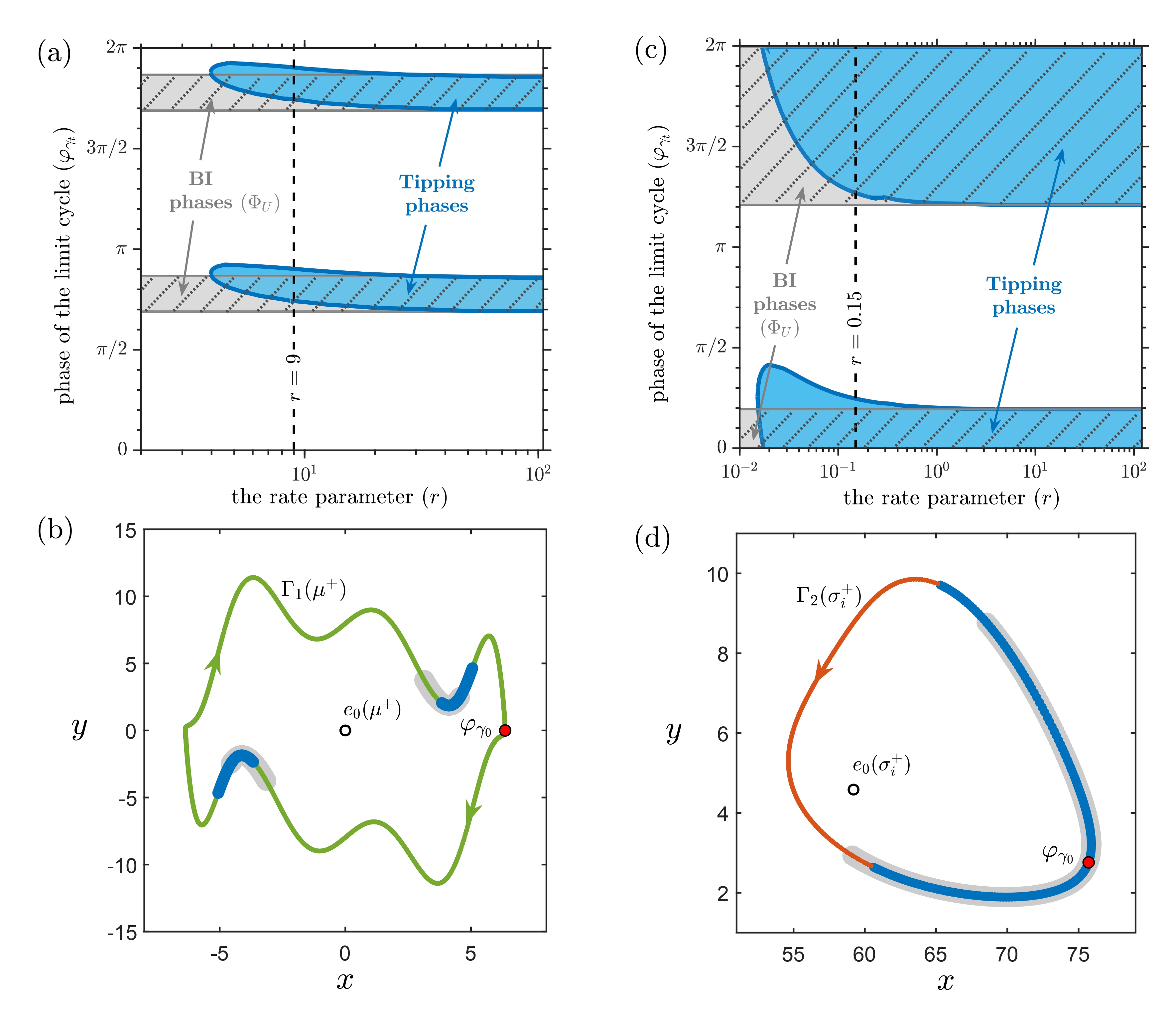}
    \caption{(a) Pace versus phase for the birhythmic van der Pol model \eqref{vdp_with_xy} with monotone time-dependent external input $\mu(rt)$ \eqref{monotone_shift}. The blue region denotes the tipping phases, and the overlapping grey-shaded region denotes basin unstable phases $(BI)$. (b) Phase portrait of large amplitude base stable limit cycle $\Gamma_1$ at $\mu^+ = 1.52$ with $BI$ phases $\Phi_U$ (grey) and tipping phases (blue) correspond to the rate value $r = 9$ (see the black dashed line in (a)). (c) Describes the same as (a) for the Decroly-Goldbeter glycolysis model \eqref{gly_model} with a monotone input $\sigma_{i}(rt)$ \eqref{monotone_shift}. (d) Small-amplitude base stable limit cycle $\Gamma_2$ at $\sigma_i^+ = 1.226$ with basin unstable $\Phi_U$ (grey) and tipping phases (blue) for the rate $r = 0.15$ (see the black dashed line in (c)).}
    \label{fig: Figure_tanh_phaserate}
\end{figure}

\begin{figure}[ht!]
    \centering
    \includegraphics[width=0.9\textwidth]{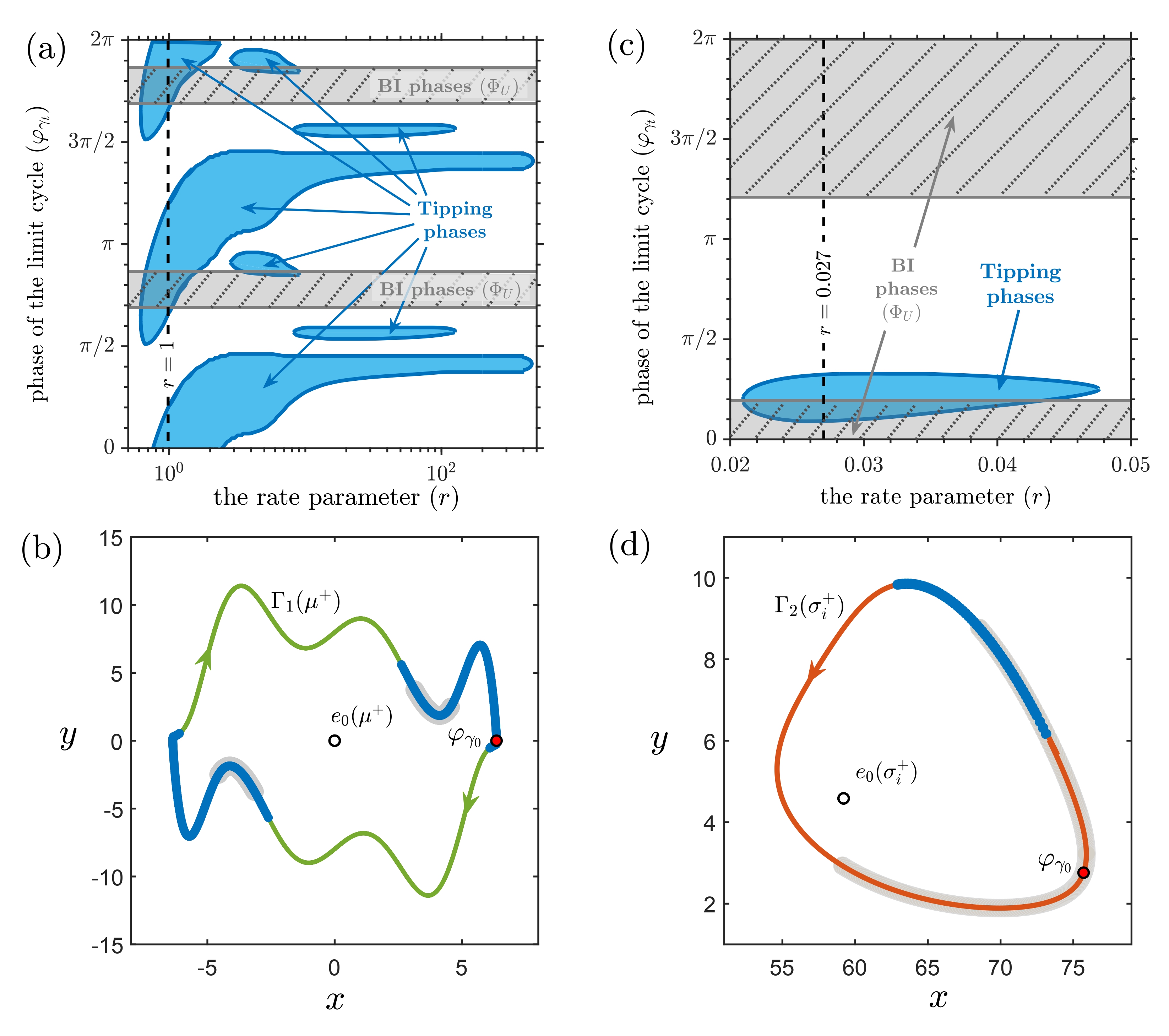}
    \caption{(a) Pace versus phase for the birhythmic van der Pol model \eqref{vdp_with_xy} with non-monotone time-dependent external input $\mu(rt)$ \eqref{nonmonotone_shift}. The blue region denotes the tipping phases, and the grey-shaded region denotes basin unstable phases $(BI)$. (b) Phase portrait of large amplitude base stable limit cycle $\Gamma_1$ at $\mu^+ = 1.52$ with $BI$ phases $\Phi_U$ (grey) and tipping phases (blue) correspond to the rate value $r = 1$ (see the black dashed curve in (a)). (c) describes the same as (a) for the Decroly-Goldbeter glycolysis model \eqref{gly_model} with a non-monotone input $\sigma_{i}(rt)$ \eqref{nonmonotone_shift}. (d) Small-amplitude base stable limit cycle $\Gamma_2$ at $\sigma_i^+ = 1.226$ with basin unstable $\Phi_U$ (grey) and tipping phases (blue) for the rate $r = 0.027$ (see the black dashed curve in (c)).}
    \label{fig: Figure_sech_phaserate}
\end{figure}

While we have discussed the concept of RP-tipping for a fixed initial condition on the base stable limit cycle, it is crucial to examine the relationship between the rate of external inputs and all initial conditions concerning the phase of the base stable limit cycles. This section explores this relationship through Figs.~\ref{fig: Figure_tanh_phaserate} and \ref{fig: Figure_sech_phaserate}. We begin by discussing this in the context of both models \eqref{vdp_with_xy} and \eqref{gly_model} with a monotone input \eqref{monotone_shift}. Fig.~\ref{fig: Figure_tanh_phaserate}(a) illustrates the relationship between the rate parameter $r$ and the phase $\varphi_{\gamma_t}$ of the stable limit cycle $\Gamma_1(\mu^+)$ in the birhythmic van der Pol model. The figure shows two distinct regions corresponding to tipping phases (blue) and basin unstable (BI) phases (grey). Due to the symmetry of the limit cycle, these regions are repeated twice. As $r \rightarrow \infty$, the monotone input steepens significantly, approaching a step function, which results in the convergence of tipping and basin unstable phases. To further illustrate this, we consider a fixed rate $r = 9$. Fig.~\ref{fig: Figure_tanh_phaserate}(b) displays the base limit cycle $\Gamma_1(\mu^+)$, highlighting two sets of phases concerning tipping (blue) and basin instability (grey).

Similarly, for higher values of the rate of an input $r$, Fig.~\ref{fig: Figure_tanh_phaserate}(c) illustrates the alignment of tipping phases (blue) with basin unstable phases (grey) of the base limit cycle $\Gamma_2(\sigma_i^+)$ in the glycolysis model. The fragmentation of the regions is due to the position of the initial phase $\varphi_{\gamma_0}$. Fig.~\ref{fig: Figure_tanh_phaserate}(d) demonstrates the phase portrait corresponding to Fig.~\ref{fig: Figure_tanh_phaserate}(c) for $r = 0.15$, showing the distribution of tipping and BI phases along the base limit cycle $\Gamma_2(\sigma_i^+)$.

For the non-monotone input \eqref{nonmonotone_shift}, tipping phases and BI phases of the base limit cycle exhibit distinct characteristics in each of the models. Unlike the monotone case, these phases do not overlap as $r \rightarrow \infty$; instead, the tipping phases terminate at a finite rate value. This occurs due to the same asymptotic limit of the non-monotone input in both directions. For a higher rate, the peak of the input undergoes a sharp transition and then quickly returns, creating the impression that there was no transition. Fig.~\ref{fig: Figure_sech_phaserate}(a) shows the fragmented regions of tipping phases (blue) and a constant strip of BI phases (grey) of the base stable limit cycle $\Gamma_1(\mu^+)$ in the birhythmic van der Pol model. 
The shaded regions repeat due to the symmetry of the limit cycle $\Gamma_1$. Fig.~\ref{fig: Figure_sech_phaserate}(b) presents the base limit cycle $\Gamma_1(\mu^+)$, illustrating two sets of phases related to tipping (blue) and basin instability (grey) for $r = 1$. Similarly, for the Decroly-Goldbeter glycolysis model, Fig.~\ref{fig: Figure_sech_phaserate}(c) depicts the tipping and BI phases. Compared to Fig.~\ref{fig: Figure_sech_phaserate}(a), the blue-shaded region here is smaller, constrained between 0.02 and 0.05 rate values. For $r = 0.027$, Fig.~\ref{fig: Figure_sech_phaserate}(d) exhibits the phase portrait corresponding to Fig.~\ref{fig: Figure_sech_phaserate}(c),   illustrating the distribution of tipping and BI phases along the base limit cycle $\Gamma_2(\sigma_i^+)$.
Intuitively, the observed discrepancy between the (blue) tipping phases and the (grey) basin unstable phases arises for non-monotone inputs for two reasons. First, RP-tipping can occur from states that are in the basin of attraction but not close enough to the limit cycle, in which case one would need to consider the basin instability of such points. Second, there are so-called rescue events~\cite{alkhayuon2023stochastic}. In such an event, the system moves to the other basin of attraction, but only temporarily, because rapid trend reversal in the input pushes the system back to the initial basin of attraction and prevents tipping.

\section{Conclusions \label{Sec: conclusion}}

Owing to the significance of evaluating a threshold in the form of a critical rate of parameters for a sudden collapse in a system state in various real-world systems, research and studies on rate-induced critical transitions are blooming over the current period. In this work, we analysed two distinct birhythmic systems of different nonlinear structures with time-varying external inputs to study RP-tipping between two stable limit cycles, crossing an unstable limit cycle that acts as a basin boundary separating the basin of attraction for both the stable limit cycles. Since the limit cycles are the primary source for the rate-induced phase-tipping, it is essential to define the phase. Due to non-uniformity in the shape of the limit cycle in the birhythmic van der Pol model, unlike the glycolysis model, we defined the phase of the limit cycle using the time instances of the system on the limit cycle. RP-tipping happens not from every phase on the limit cycle but from certain phases depending on the rate $r$ of external input. We considered two types of parameter shifts - monotone and non-monotone - of the chosen parameters in both models to vary over time. We numerically investigated the existence of RP-tipping using two-dimensional tipping diagrams with the parameters in the shift - magnitude and rate - fixing the value of the parameter $t_c$. To make our results robust, we used the concept of partial basin instability of the limit cycle $\Gamma(p)$ on the fixed-in-time input parameters along a parameter path $\Delta_p$, which is a set-theoretic property of the autonomous frozen birhythmic systems. It is proposed that if $\Gamma(p)$ is partially basin unstable on a parameter path $\Delta_p$, then there exists a parameter setting on $\Delta_p$ which assures RP-tipping crossing the critical rate. We presented the partial basin instability of $\Gamma_1$ in the birhythmic van der Pol model and $\Gamma_2$ in the glycolysis for our investigation on RP-tipping.

The most challenging component of partial basin instability of $\Gamma(p)$ and RP-tipping is determining the suitable parameter path $\Delta_p$ on which the shifts evolve with respect to time. It is not necessary to have partial basin instability of $\Gamma(p)$ with the randomly chosen parameter paths. Choosing the direction of drift of the parameters on $\Delta_p$ is also as important as choosing the parameter path. In our case, for both models, we decreased the parameters from $p^+$ to $p^-$ to witness the basin instability. To distinguish the parameter paths on the birhythmic region of two-parameter bifurcation diagrams, for the birhythmic van der Pol model, the path is chosen vertically, parallel to $\mu$-axis with fixed $d$ whereas, for the glycolysis model, the parameter path is chosen diagonally by changing both the substrate input $v$ and the maximal rate of recycling of the product $\sigma_i$ in a straight line. However, one may determine other parameter paths. Further, in the birhythmic van der Pol model for a longer parameter path and a fixed rate, we illustrated a series of RP-tipping between limit cycles using a distinct non-monotone input in the form of impulse. At last, the relationship between tipping phases and basin unstable phases is studied using the pace versus phase diagrams for the forcings in both models.

Taking into account various real-life scenarios, it is becoming apparent that R-tipping is just as crucial as B-tipping. It is of the utmost importance to develop alternative techniques, like early warning signals for different B-tippings, to forecast the occurrence of R-tipping, and then RP-tipping. It is also important to verify R-tipping in models with diverse potential dynamics \citep{ashwin2021physical}. Considering nonlinear models with the coexistence of chaotic attractors for R-tipping will also be quite interesting to study. This present analysis can be made more general by developing the concepts of basin instability and determining the phases for the chaotic attractors. 

\section*{Data accessibility.}
The article has no additional data. \\
Supplementary material is available online:\\ (\url{https://github.com/RaviKumarK97/Electronic_supplementary_RP-tipping}).

\section*{Declaration of AI use.} We have not used AI-assisted technologies in creating this article.

\section*{Authors' contributions.} All authors contributed equally.\\
All authors gave final approval for publication and agreed to be held accountable for the work performed therein.

\section*{Conflict of interest declaration.} We declare we have no competing interests.

\section*{Funding.} Ravi Kumar K has received funding from IIT Ropar.

\section*{Acknowledgments.}
P.S.D. acknowledges financial support from the Science \& Engineering Research Board (SERB), Government of India (Grant number: CRG/2022/002788). 


\end{document}